\documentclass[11pt]{article}

\usepackage{multirow}
\usepackage{enumitem}
\usepackage{caption}
\usepackage[colorlinks, linkcolor=blue, anchorcolor=red, citecolor=blue, CJKbookmarks=true]{hyperref}
\usepackage{epsfig, epstopdf, graphicx, floatrow, rotating, amssymb, bbm, amsthm, amsbsy, amsmath, amsfonts, mathrsfs, dsfont, makeidx, color}
\usepackage{bm}
\usepackage{algorithmicx}
\usepackage{algpseudocode}
\usepackage{natbib}
\usepackage{graphicx}
\usepackage{float}
\usepackage{subfigure}
\usepackage{xr}
\usepackage[strict]{changepage}

\usepackage[left=4cm,right=2.5cm,top=2cm,bottom=4cm,includehead,includefoot,headheight=15pt]{geometry}
\allowdisplaybreaks

\usepackage{setspace}
\graphicspath{{images/}}
\numberwithin{equation}{section}

\theoremstyle{plain}
\newtheorem{theorem}{Theorem}
\newtheorem{corollary}{Corollary}[theorem]
\newtheorem{proposition}{Proposition}
\newtheorem{lemma}{Lemma}[section]
\newtheorem{claim}{Claim}[section]
\newtheorem{definition}{Definition}
\newtheorem{assumption}{Assumption}[section]

\theoremstyle{remark}
\newtheorem*{remark}{Remark}
\newtheorem{example}{Example}[section]

\DeclareMathOperator*{\argmin}{argmin}

\DeclareMathOperator*{\var}{Var}
\DeclareMathOperator*{\cov}{Cov}
\DeclareMathOperator*{\lip}{Lip}

\newcommand\norm[1]{\left\lVert#1\right\rVert}

\newcommand\bracketL[1]{\left\{#1\right\}}

\newcommand{\expct}{\mathbb{E}}
\newcommand{\prob}{\mathbb{P}}
\newcommand{\cA}{\mathcal{A}}
\newcommand{\cB}{\mathcal{B}}

\newcommand{\cF}{\mathcal{F}}
\newcommand{\cG}{\mathcal{G}}

\newcommand{\cJ}{\mathcal{J}}

\newcommand{\cN}{\mathcal{N}}
\newcommand{\cX}{\mathcal{X}}

\newcommand{\fX}{\mathfrak{X}}
\newcommand{\fY}{\mathfrak{Y}}

\newcommand{\bL}{\mathbb{L}}
\newcommand{\bN}{\mathbb{N}}
\newcommand{\bR}{\mathbb{R}}

\newcommand{\bT}{\mathbb{T}}

\newcommand{\bZ}{\mathbb{Z}}
\newcommand{\img}{\mathbf{i}}

\newcommand{\convp}{\overset{p}{\to}}
\newcommand{\convd}{\overset{d}{\to}}
\newcommand{\convLone}{\overset{L^1}{\to}}

\newcommand{\bigO}{\mathcal{O}}

\begin{document}

\title{Asymptotic theory and inference for non-stationary and non-mixing triangular arrays of random fields}

\author{
Yue Pan$^{1,2}$ \quad Jiazhu Pan$^{2}$ \\
\textsl{\footnotesize $^{1}$School of Economics, Anhui University, Hefei 230601, China} \\
\textsl{\footnotesize $^{2}$Department of Mathematics and Statistics, University of Strathclyde, Glasgow G1 1XH, UK}
}
\date{}
\maketitle

\begin{abstract}
We develop a unified asymptotic framework for non-stationary, non-mixing triangular arrays of random fields on multi-dimensional lattices under row-uniform $\eta$-weak dependence. We establish the preservation of weak dependence under locally Lipschitz transformations and spatially and row-wise heterogeneous Bernoulli shifts, obtaining explicit dependence bounds determined by innovation dependence and functional sensitivity. Under uniform moment conditions, polynomially decaying $\eta$-dependence coefficients, and a nondegenerate aggregate variance condition, we prove a law of large numbers and a central limit theorem. The scope of the framework is illustrated through several examples that highlight its ability to accommodate non-stationarity without requiring mixing assumptions. As an application, we develop parameter inference procedures for a spatio-temporal model with a dynamic network structure, thereby demonstrating the practical relevance of the proposed asymptotic theory.
\end{abstract}


%\textbf{JEL Classification}: C10, C31, C33.


\section{Introduction}\label{section_introduction}

Spatio-temporal models have drawn much attention in science, social science, and econometrics because of their merits in capturing dependence over both time and space. However, samples from these models often form multi-dimensional panels whose finite sampling regions expand over space and time, while their distributions and dependence structures may also change with the sample size. The asymptotic theory used in the inference of univariate or fixed-dimensional time series therefore does not directly apply. A spatio-temporal model can be regarded as a random process running on a multi-dimensional lattice, i.e., a random field. Thus, a non-stationary expanding sample can be represented by a triangular array of random fields. We establish asymptotic theory for such arrays that provide useful tools in the inference of spatio-temporal models.

Asymptotic theory for random fields has been extensively studied in the literature. \cite{jenish2009} proposed asymptotic theory for triangular arrays of random fields under $\alpha$- and $\phi$-mixing. Compared to previous asymptotic theory for mixing random fields \citep{bolthausen1982, guyon1995, dedecker1998}, their asymptotic theory is more general in the sense that it accommodates non-stationarity. In this paper, we show the limitation of mixing through two examples: Section \ref{subsection_non-stationary_bernoulli_array} provides a Bernoulli autoregressive triangular array of random fields that is non-stationary both within each row and across rows. Section \ref{subsection_dynamic_network_bernoulli_array} allows innovations to propagate through a sequence of local network matchings that changes over time. In both examples, the corresponding $\alpha$- and $\phi$-mixing coefficients remain bounded away from zero at every distance. These examples identify a class of spatio-temporal models that lies beyond the scope of existing asymptotic theory based on mixing conditions.

To address this restriction of mixing, \cite{doukhan1999} introduced the concept of weak dependence, which can be extended to random fields. See Chapter 2 in \cite{dedecker2007} for details. However, existing asymptotic theory for weakly dependent random fields either requires stationarity \citep{dedecker2007, machkouri2013, curato2022} or is only established for single-indexed sequences \citep{neumann2013, merlevede2019}. For example, \cite{machkouri2013} proposed a CLT for stationary random fields that are Bernoulli shifts of IID innovations, which provides an important stationary benchmark under a different dependence condition. \cite{neumann2013} and \cite{merlevede2019} proposed a CLT and a functional CLT, respectively, for non-stationary triangular arrays of random variables rather than triangular arrays of random fields with two or more indices. Since they are both limited to random sequences along a single time index, their results are not applicable to the triangular arrays of random fields with multi-dimensional indices that we discuss later. In Theorem \ref{theorem_lln} and Theorem \ref{theorem_clt}, we establish a law of large numbers and a central limit theorem directly under $\eta$-weak dependence for non-stationary and non-mixing triangular arrays of random fields. The moment conditions and dependence rate conditions are imposed uniformly over the rows of the array. Corollary \ref{corollary_clt} further extends the CLT to triangular arrays of vector-valued random fields, which is useful for scores and estimating equations.

To apply the established asymptotic theory to the inference of dynamic spatio-temporal models, it is important to establish how weak dependence is preserved under nonlinear transformations and infinite shifts. Such hereditary properties have been proved under certain conditions. See, for example, Proposition 2.4 in \cite{curato2022}. In Proposition \ref{proposition_transformation}, we extend this property to locally Lipschitz transformations of triangular arrays of random fields. In Proposition \ref{proposition_BS}, we also show that weak dependence can be preserved under infinite shifts whose functional form and sensitivity coefficients may vary with both location and the array row. This extends the stationary heredity results in \cite{dedecker2007}, \cite{doukhanWintenberger2007}, and \cite{doukhan2007} to non-stationary triangular arrays of random fields. Section \ref{subsection_non-stationary_ambit_mma} applies this result to a L\'{e}vy-driven triangular array of random fields with non-stationary ambit and mixed moving-average representations. Both the kernel loading and the causal ambit geometry vary over locations and rows, while the dependence coefficient decays exponentially up to a polynomial factor.

Furthermore, we use the dynamic-network model in Section \ref{subsection_dynamic_network_bernoulli_array} to illustrate how the heredity results and asymptotic theory lead to parameter inference. Theorem \ref{theorem_dynamic_network_estimator} establishes the consistency and asymptotic normality of the least squares estimator and the consistency of the sandwich covariance estimator. With these results, we build a direct link from weak dependence in a non-stationary and non-mixing spatio-temporal process with a dynamic network structure to feasible parameter inference. The same argument also provides the dependence and limit-theory ingredients for smooth M-, QML-, and GMM procedures with controlled local propagation, including models related to network vector autoregressions \citep{zhu2017}, nonlinear network autoregressions \citep{armillotta2023}, and time series on dynamic networks \citep{krampe2019}.

The rest of this paper is organized as follows. In Section \ref{section_weak_dependence}, we introduce the concept of $\eta$-weak dependence for triangular arrays of random fields, provide the non-stationary and non-mixing examples, and investigate its heredity under transformations and infinite shifts. Our LLN and CLT for weakly dependent triangular arrays of random fields are presented in Section \ref{section_limit_theorem}. In Section \ref{section_dynamic_network_estimation}, we establish the asymptotic properties of nonlinear least-squares estimation for a spatio-temporal model with a dynamic network structure. Section \ref{section_conclusion} concludes the paper. Proofs and verifications of the main results are collected in Appendix \ref{appendix_proofs}.

We use the following notation throughout. Let $\bN=\{1,2,\ldots\}$ and $\bN_0=\{0,1,2,\ldots\}$. The symbols $\expct$ and $\prob$ denote expectation and probability, respectively, and $\mathbf 1\{A\}$ denotes the indicator of an event or condition $A$. Unless otherwise specified, $\norm{\cdot}$ denotes a norm, $\norm{X}_p=(\expct|X|^p)^{1/p}$ denotes the $\bL^p$ norm of a random variable, and $|A|_c$ denotes the cardinality of a finite set $A$. We write $\lfloor x\rfloor$ and $\lceil x\rceil$ for the floor and ceiling of $x$, and $\mathcal N(\mu,\Sigma)$ for a Gaussian distribution with mean $\mu$ and covariance matrix $\Sigma$. Throughout, $C\in(0,\infty)$ denotes a generic constant that may change from line to line, whereas fixed constants are distinguished by descriptive subscripts or superscripts. When a generic constant in $(0,1)$ is needed, it is denoted by $\zeta$.

\section{Weakly dependent triangular arrays of random fields with non-stationarity}\label{section_weak_dependence}

Random fields are random processes running on multi-dimensional lattices. Considering a metric space $(\bT, \rho)$, one could easily define on $\bT$ an infinitely countable lattice $D \subset \bT$, which satisfies the following assumption throughout this section:
\begin{assumption}\label{assumption_minimum_distance}
Defined on the metric space $(\bT, \rho)$, the lattice $D \subset \bT$ is infinitely countable. For any $i, j \in D$ such that $i\neq j$, it is assumed that $\rho(i,j) \geq 1$.
\end{assumption}
\noindent This minimum distance assumption ensures that the sample size can grow as the finite sampling regions expand in $D$. A simple example satisfying Assumption \ref{assumption_minimum_distance} is the $d$-dimensional Euclidean space $\bT=\bR^d$ with the infinite lattice $D=\bZ^d$, whose minimum distance is 1. Let $(D_n)_{n\geq1}$ be a sequence of finite sampling regions $D_n\subset D$. For each $n$, $\{X_{i,n}:i\in D_n\}$ is a random field taking values in a Banach space $(\cX,\|\cdot\|)$; the double-indexed collection $\{X_{i,n}:i\in D_n,n\geq1\}$ is a triangular array of random fields.

In their Definition 2.2, \cite{dedecker2007} defined the $(\cF, \cG, \Psi)$-coefficients that measure the dependence between two separated groups of random variables on $\bZ$, and remarked that the definition extends to general metric index sets. Random-field versions of the $\eta$-coefficient were subsequently used by \cite{doukhan2007} and \cite{curato2022}. In Section \ref{subsection_weak_dependence}, we adopt this standard coefficient and formulate it uniformly over the rows of a non-stationary triangular array of random fields indexed by any lattice $D$ satisfying Assumption \ref{assumption_minimum_distance}.

\subsection{Weak dependence}\label{subsection_weak_dependence}

Let $\cF_u$ and $\cG_v$ denote two classes of functions from $\cX^u$ to $\bR$ and $\cX^v$ to $\bR$ respectively. If $\Psi$ is some function mapping from $\cF_u \times \cG_v$ to $\bR_+$, and $\fX$ and $\fY$ are two arbitrary random variables in $\cX^u$ and $\cX^v$, then we can define the measurement of dependence between $\fX$ and $\fY$ by
\begin{equation}\label{eta}
    \eta(\fX, \fY) = \sup\left\{\frac{|\cov(f(\fX), g(\fY))|}{\Psi(f,g)}: f\in\cF_u, g\in\cG_v\right\},
\end{equation} where $\Psi(f, g) = u\|g\|_{\infty}\lip(f) + v\|f\|_{\infty}\lip(g)$.

Given any $U_n\subset D_n$ with cardinality $|U_n|_c = u$ and $V_n\subset D_n$ with cardinality $|V_n|_c = v$, let $\fX_{U_n} = (X_{i,n})_{i \in U_n}$, $\fX_{V_n} = (X_{i,n})_{i \in V_n}$. Then the dependence coefficient of the triangular array of random fields $\{X_{i,n}: i\in D_n, n\geq 1\}$ is defined by
\begin{equation}\label{dependence_coef}
    \eta_{n,u,v}(r) = \sup\left\{\eta(\fX_{U_n}, \fX_{V_n}): |U_n|_c = u, |V_n|_c = v, \rho(U_n, V_n) \geq r\right\},
\end{equation} where $\rho(U_n, V_n) := \min_{i \in U_n, j \in V_n} \rho(i, j)$ measures the distance between $U_n$ and $V_n$.
\begin{remark}
For any functions $f \in \cF_u$ and $g \in \cG_v$ such that $\Psi(f,g) < \infty$, the inequality below follows directly from \eqref{dependence_coef}:
\begin{equation}\label{cov_ineq_bounded}
    |\cov(f(\fX_{U_n}), g(\fX_{V_n}))| \leq C \eta_{n,u,v}(\rho(U_n,V_n))
\end{equation} for some constant $0 < C < \infty$ that is related to $f$, $g$ and $\Psi$.
\end{remark}
\begin{remark}
Similar to Definition 1 in \cite{jenish2009}, we introduce the following notation: \[\bar\eta_{u,v}(r) = \sup_{n}\eta_{n,u,v}(r), \qquad \bar\eta(r) = \sup_{u,v}\bar\eta_{u,v}(r).\]
\end{remark}

Now we are ready to give a formal definition of weak dependence as follows:
\begin{definition}\label{def_of_weak_dep}
The triangular array of random fields $\{X_{i,n}: i\in D_n, n\geq 1\}$ in Banach space $(\cX, \|\cdot\|)$ is said to be $\eta$-weakly dependent if $\lim_{r\to\infty}\bar\eta(r) = 0$.
\end{definition}

\eqref{dependence_coef} is a standard form of $\eta$-weak-dependence coefficients. The original formulation in \cite{doukhan1999} compares time-ordered blocks of a one-dimensional sequence, and much of the associated asymptotic theory is developed under stationarity. Existing random-field formulations and asymptotic theory, including \cite{doukhan2007} and \cite{curato2022}, likewise focus on stationary fields. Conversely, \cite{neumann2013} allows possibly non-stationary triangular arrays under weak dependence, but for one-dimensional sequences. Thus, the role of our formulation is to combine non-stationarity, a triangular-array structure, and multi-dimensional indexing in one row-uniform condition.

\subsection{A non-stationary Bernoulli autoregressive triangular array of random fields}\label{subsection_non-stationary_bernoulli_array}

\cite{jenish2009} already accommodate non-stationary triangular arrays of random fields on multi-dimensional, possibly irregular lattices. Their Definition 1 also takes a supremum over the row index, but their CLT requires uniform decay of $\alpha$- or $\phi$-mixing coefficients defined through the $\sigma$-fields generated by separated blocks. Since mixing controls all measurable events, it may fail when a distant observation retains an exactly recoverable discrete innovation. In contrast, $\eta$-weak dependence controls bounded Lipschitz transformations, for which the numerical effect of that innovation may still decay rapidly. The construction below makes this distinction explicit. Its binary autoregressive mechanism is classical \citep{andrews1984,doukhan1999}; the purpose of the triangular-array construction is to illustrate the scope of Definition \ref{def_of_weak_dep}.

Let $D=\bZ^2$ be equipped with
\[
\rho((i,t),(j,s))=\max\{|i-j|,|t-s|\},
\qquad
D_n=\{1,\ldots,4n\}^2.
\]
For $(i,t)\in D$ and $n\geq1$, let the variables $\varepsilon_{i,t,n}$ be mutually independent and
\begin{equation}\label{example_bernoulli_prob}
\varepsilon_{i,t,n}\sim\operatorname{Bernoulli}(p_{i,t,n}),
\qquad
p_{i,t,n}=\frac12+\frac14
\sin\left(\frac{\pi i}{2n}\right)
\cos\left(\frac{\pi t}{2n}\right).
\end{equation}
In particular, $1/4\leq p_{i,t,n}\leq3/4$. Consider the Bernoulli autoregression
\[
    X_{i,t,n}
    =\frac12X_{i,t-1,n}
    +\frac12\varepsilon_{i,t,n}.
\]
Then the following statements hold.
\begin{enumerate}
    \item[(i)] The autoregression has the unique bounded causal solution
    \begin{equation}\label{example_bernoulli_array}
        X_{i,t,n}
        =\sum_{k=0}^{\infty}2^{-(k+1)}
        \varepsilon_{i,t-k,n},
        \qquad
        0\leq X_{i,t,n}\leq1.
    \end{equation}
    The triangular array of random fields $\bracketL{X_{i,t,n}: (i,t)\in D_n, n\geq1}$ is non-stationary within each row, and its law changes with $n$.
    \item[(ii)] It is $\eta$-weakly dependent in the sense of Definition \ref{def_of_weak_dep}, with
    \begin{equation}\label{example_eta_bound}
        \eta_{n,u,v}(r)\leq 2^{1-\lceil r\rceil}
        \quad\text{and hence}\quad
        \bar\eta(r)\leq2^{1-\lceil r\rceil},
        \qquad r\geq1.
    \end{equation}
    The bound is uniform in $n,u,v$ and applies to arbitrary sets $U_n,V_n$ separated under the random-field metric.
    \item[(iii)] Let $\bar{\alpha}_{k,l}(r)$ and $\bar{\phi}_{k,l}(r)$ denote the uniform-in-$n$ triangular-array coefficients in Definition 1 of \cite{jenish2009}. For every integer $r\geq1$,
    \begin{equation}\label{example_mixing_lower_bound}
        \bar{\alpha}_{1,1}(r)\geq\frac{3}{16},
        \qquad
        \bar{\phi}_{1,1}(r)\geq\frac14.
    \end{equation}
    Consequently, the array is neither $\alpha$- nor $\phi$-mixing.
\end{enumerate}
The verification of above claims is given in Section \ref{proof_example_non-stationary_bernoulli_array}.

\subsection{A spatio-temporal model with a dynamic network structure}\label{subsection_dynamic_network_bernoulli_array}

In the preceding example, $X_{i,t,n}$ is affected by its lag-1 history $X_{i,t-1,n}$. We next show that the distinction between weak dependence and mixing persists when $X_{i,t,n}$ is also affected by its spatial neighbors through a dynamic network.

Let $D=\bZ^2$ have the same maximum metric as in Section \ref{subsection_non-stationary_bernoulli_array}, and let $D_n=\{1,\ldots,4n\}^2$. For each node $i\in\bZ$, let
\begin{equation}\label{example_dynamic_network_matchings}
\begin{aligned}
\pi^{(0)}(i)&=
\begin{cases}
i+1,&i\text{ is even},\\
i-1,&i\text{ is odd},
\end{cases}
&
\pi^{(1)}(i)&=
\begin{cases}
i-1,&i\text{ is even},\\
i+1,&i\text{ is odd}.
\end{cases}
\end{aligned}
\end{equation}
Let $\Gamma\subset\bR$ be compact. Given a deterministic sequence $(q_t)_{t\in\bZ}$ and a threshold $\gamma\in\Gamma$, for two arbitrary nodes $i,j\in\bZ$, set
\begin{equation}\label{example_dynamic_network_regime}
    \pi_t^{\gamma}=\pi^{(\mathbf 1\{q_t>\gamma\})},
    \qquad
    w_{ij,t}(\gamma)=\mathbf 1\{j=\pi_t^{\gamma}(i)\}.
\end{equation}
Nodes $i$ and $j$ are connected if $w_{ij,t}(\gamma) = 1$. Thus the network is dynamic since the relations between nodes change whenever $q_t$ crosses $\gamma$.

Let $d_\beta\in\bN$, and let $\cB\subset\bR^{d_\beta}$ be compact. For each $n$ and $\beta\in\bR^{d_\beta}$, let the variables $\varepsilon_{i,t,n}$ be mutually independent over $(i,t)\in D$ and satisfy
\begin{equation}\label{example_dynamic_network_prob}
    \varepsilon_{i,t,n}
    \sim\operatorname{Bernoulli}(p_{i,t,n}(\beta)),
    \qquad
    0<p_0\leq p_{i,t,n}(\beta)\leq1-p_0<1
\end{equation}
uniformly in $i,t,n$ and $\beta$.
We specialize the innovation probability to the logistic specification
\begin{equation}\label{estimation_dynamic_network_probability}
    p_{i,t,n}(\beta)
    =\frac{\exp(z_{i,t,n}'\beta)}{1+\exp(z_{i,t,n}'\beta)}
\end{equation}
where the deterministic covariates $z_{i,t,n}\in\bR^{d_\beta}$ are uniformly bounded. This permits the law of innovation $\varepsilon_{i,t,n}$ to vary with the $i$, $t$, and $n$.

For the parameter
\[
    \vartheta=(a,\beta,\gamma)\in
    \Theta=[\underline a,\bar a]\times\cB\times\Gamma,
    \qquad 0<\underline a\leq\bar a<\frac12,
\] consider the dynamic-network autoregression
\begin{equation}\label{example_dynamic_network_model}
\begin{aligned}
X_{i,t,n}
&=a\sum_{j\in\bZ}w_{ij,t}(\gamma)X_{j,t-1,n}
  +\varepsilon_{i,t,n}.
\end{aligned}
\end{equation}
Define the backward network path by
\begin{equation}\label{example_dynamic_network_ancestor}
    A_0^\gamma(i,t)=i,
    \qquad
    A_{k+1}^\gamma(i,t)
    =\pi_{t-k}^\gamma\bigl(A_k^\gamma(i,t)\bigr),
    \qquad k\geq0.
\end{equation}
Thus $A_k^\gamma(i,t)$ is the node reached by tracing the network
$k$ periods backward from $(i,t)$; in particular,
$A_1^\gamma(i,t)=\pi_t^\gamma(i)$.

The following statements hold.
\begin{enumerate}
    \item[(i)] Equation \eqref{example_dynamic_network_model} has the unique bounded causal solution
    \begin{equation}\label{example_dynamic_network_causal_solution}
        X_{i,t,n}
        =\sum_{k=0}^{\infty}a^k
        \varepsilon_{A_k^\gamma(i,t),t-k,n},
        \qquad
        0\leq X_{i,t,n}\leq\frac{1}{1-a}.
    \end{equation}
    \item[(ii)] For each $\vartheta$, let $\eta_{n,u,v}^{\vartheta}(r)$ and $\bar\eta^{\vartheta}(r)$ denote the coefficients in \eqref{dependence_coef} and Definition \ref{def_of_weak_dep}. Uniformly over the parameter space,
    \begin{equation}\label{example_dynamic_network_eta_bound}
        \sup_{\vartheta\in\Theta}\sup_n\sup_{u,v}
        \eta_{n,u,v}^{\vartheta}(r)
        =\sup_{\vartheta\in\Theta}\bar\eta^{\vartheta}(r)
        \leq\frac{2\bar a^{\lceil r/2\rceil}}{1-\bar a},
        \qquad r\geq1.
    \end{equation}
    \item[(iii)] Let $\bar\alpha_{1,1}^{\vartheta}(r)$ and $\bar\phi_{1,1}^{\vartheta}(r)$ be the corresponding uniform-in-$n$ singleton mixing coefficients of \cite{jenish2009}. For every integer $r\geq1$,
    \begin{equation}\label{example_dynamic_network_mixing_lower_bound}
        \inf_{\vartheta\in\Theta}\bar\alpha_{1,1}^{\vartheta}(r)
        \geq p_0(1-p_0),
        \qquad
        \inf_{\vartheta\in\Theta}\bar\phi_{1,1}^{\vartheta}(r)
        \geq p_0.
    \end{equation}
\end{enumerate}
The verification of these claims is given in Section \ref{proof_example_dynamic_network_bernoulli_array}.

This example therefore demonstrates that Definition \ref{def_of_weak_dep} accommodates non-stationary triangular arrays of random fields generated through parameterized dynamic networks even when the corresponding mixing coefficients do not decay.

\subsection{Heredity of weak dependence}\label{section_heredity}

In this section, we will show that the $\eta$ weak dependence can be preserved under locally Lipschitz transformations and infinite shifts.

Proposition \ref{proposition_transformation} below is a natural extension of Proposition 2.1 and Proposition 2.2 in \cite{dedecker2007} to triangular arrays of random fields. It shows that weak dependence is inherited under transformations satisfying condition \eqref{lipschitz_condition_local}. A simple example is any Lipschitz-continuous function when $a = 1$.
\begin{proposition}\label{proposition_transformation}
Let $\{X_{i,n}: i\in D_n, n\geq 1\}$ be a triangular array of $\bR^{d_x}$-valued random fields with finite $p$-th moment for some $p > 1$, and let $H: \bR^{d_x} \mapsto \bR$ be a function such that
\begin{equation}\label{lipschitz_condition_local}
    |H(x) - H(y)| \leq c \|x - y\|(\|x\|^{a-1} + \|y\|^{a-1})
\end{equation} for some $c \in (0,+\infty)$, $a \in [1, p)$, and any $x, y \in \bR^{d_x}$. Let the triangular array of random fields $\{Y_{i,n}: i\in D_n, n\geq 1\}$ be defined by $Y_{i,n} = H(X_{i,n})$. If the triangular array of random fields $\{X_{i,n}: i\in D_n, n\geq 1\}$ is weakly dependent with coefficients $\bar\eta^X(r)$, then the triangular array of random fields $\{Y_{i,n}: i\in D_n, n\geq 1\}$ is also weakly dependent with $\bar\eta^Y(r) \leq C \bar\eta^X(r)^{\frac{p-a}{p-1}}$ for some constant $C > 0$.
\end{proposition}

For the heredity of weak dependence under shifts, we consider $D\subset\bZ^d$, equipped with distance measure $\rho(i,j)$ for any $i,j\in D$. Let $\{\varepsilon_i: i\in D\}$ be a $\bR$-valued random field on $D$. Let $H_{i,n}: \bR^{D} \mapsto \bR$ be a measurable function, and define the triangular array of random fields $\{X_{i,n}: i\in D_n, n\geq 1\}$ by $X_{i,n} := H_{i,n}((\varepsilon_j)_{j\in D})$. For each $h \in \bN_0$, and for any $(x_j)_{j\in D}$ and $(y_j)_{j\in D}$ such that $x_j \neq y_j$ if and only if $\rho(i,j) = h$, $H_{i,n}$ satisfies that
\begin{equation}\label{lipschitz_condition_BS}
\begin{aligned}
    &|H_{i,n}((x_j)_{j\in D}) - H_{i,n}((y_j)_{j\in D})|\\
    \leq &B_{i,n}(h)(\max_{\rho(i,j)\neq h}|x_j|^l \vee 1)\sum_{\rho(i,j)=h}|x_j-y_j|
\end{aligned}
\end{equation} almost surely, where $l \geq 0$ and $\{B_{i,n}(h): i\in D_n, n\geq 1\}$ are positive constants satisfying that
\begin{equation}\label{coefficient_condition_BS}
    C_B := \sup_{n\geq 1}\sup_{i\in D_n}\sum_{h=0}^{\infty}B_{i,n}(h)h^{d-1} < \infty.
\end{equation} In Proposition \ref{proposition_BS} below, we investigate the preservation of weak dependence from the $\{\varepsilon_i: i\in D\}$ to $\{X_{i,n}: i\in D_n, n\geq 1\}$.
\begin{proposition}\label{proposition_BS}
    Let $\{X_{i,n}\in\bR: i\in D_n, n\geq 1\}$ be a triangular array of random fields generated by Bernoulli shifts of $\{\varepsilon_i\in\bR: i\in D\}$, such that $X_{i,n} = H_{i,n}((\varepsilon_j)_{j\in D})$ and $H_{i,n}: \bR^{D} \mapsto \bR$ satisfies conditions \eqref{lipschitz_condition_BS} and \eqref{coefficient_condition_BS}. Assume that $\sup_{i\in D}\expct|\varepsilon_i|^p < \infty$ with $p > l + 1$. If the random field $\{\varepsilon_i\in\bR: i\in D\}$ is weakly dependent with coefficients $\bar\eta^\varepsilon(r)$, then the triangular array of random fields $\{X_{i,n}\in\bR: i\in D_n, n\geq 1\}$ is also weakly dependent with coefficients
    \begin{equation}\label{eq_eta_output}
        \bar\eta(r) = C \inf_{0 < s \leq \lfloor r/2\rfloor}\left\{R_B(s)\vee\left[s^d\bar\eta^\varepsilon(r-2s)^{\frac{p-1-l}{p-1}}\right]\right\},
    \end{equation} where $R_B(s) = \sup_{n\geq 1}\sup_{i\in D_n}\sum_{h\geq s}B_{i,n}(h)h^{d-1}$ is the sensitivity tail.
\end{proposition}

It is not easy to find the exact infinum in \eqref{eq_eta_output}. However, the dependence coefficients of the outputs have upper bounds in explicit forms, if the dependence coefficients and $B_{i,n}$ decay in a regular manner.
\begin{example}\label{example_wd_bound_power}
Let the dependence coefficients of the input $\bar\eta^\varepsilon(r) = \bigO(r^{-\mu})$ for some $\mu > \frac{p-1}{p-1-l}d$, and $B_{i,n}(h) = \bigO(h^{-b})$ for some $b \geq \frac{p-1-l}{p-1}\mu$. Then the dependence coefficients of the output are bounded by:
\begin{equation}\label{eq_wd_bound_power}
    \bar\eta(r) \leq C r^{d-\frac{p-1-l}{p-1}\mu}.
\end{equation}
\end{example}
\begin{example}\label{example_wd_bound_exp}
Assume that $d = 2$, let the dependence coefficients of the input $\bar\eta^\varepsilon(r) = \bigO(r^{-\mu})$ for some $\mu > 0$, and $B_{i,n}(h) = \bigO(e^{-bh})$ for some $b \geq \frac{p-1-l}{p-1}\mu$. Then the dependence coefficients of the output are bounded by:
\begin{equation}\label{eq_wd_bound_exp}
    \bar\eta(r) \leq C(\log r)^2 r^{-\frac{p-1-l}{p-1}\mu}.
\end{equation}
\end{example} \noindent The proofs of \eqref{eq_wd_bound_power} and \eqref{eq_wd_bound_exp} are given in Section \ref{eq_wd_bound_proof}.

The heredity principle for Bernoulli shifts is well established in stationary settings. For stationary sequences, \cite{doukhanWintenberger2007} and Chapter 3 of \cite{dedecker2007} show that weak dependence is preserved under infinite-memory shift functionals, including locally Lipschitz functionals of weakly dependent innovations; \cite{doukhan2007} obtains related inheritance results for stationary random fields generated by homogeneous fixed-point representations.

The main advance of Proposition \ref{proposition_BS} is to make this inheritance theory applicable to non-stationary triangular arrays of random fields: the shift functional $H_{i,n}$ and the sensitivity coefficients $B_{i,n}(h)$ may vary jointly with the spatial location $i$ and the array row $n$. This accommodation of non-stationarity, triangular-array asymptotics, multi-dimensional indexing, and spatially and row-wise heterogeneous shift mechanisms is not covered by the stationary results above and is the distinctive contribution of Proposition \ref{proposition_BS}. The next subsection places this contribution in a L\'{e}vy-driven setting by constructing a non-stationary triangular array of random fields with ambit and mixed moving-average representations that extends the stationary kernel and ambit specifications studied by \cite{curato2022}.

\subsection{A non-stationary triangular array of random fields with ambit and mixed moving-average representations}\label{subsection_non-stationary_ambit_mma}
Let $d=d_s+1$ with $d_s\geq1$, write $i=(t,x)\in\bZ\times\bZ^{d_s}$, equip $D=\bZ^d$ with the maximum metric, and take $D_n=\{1,\ldots,n\}^d$. For $j=(j_0,\ldots,j_{d_s})\in\bZ^d$, let
\[
    Q_j=(j_0-1,j_0]\times
    \prod_{k=1}^{d_s}[j_k,j_k+1),
\] be a $d$-dimensional unit spatio-temporal cell.
Let $\Lambda$ be a centered homogeneous factorisable real-valued L\'{e}vy basis on $\mathcal M\times\bR^d$ with control measure $\pi(dA)\,ds$. Choose a deterministic mixing function $q$ such that the innovation
\begin{equation}\label{example_ambit_innovations}
    \varepsilon_j
    =\int_{\mathcal M}\int_{Q_j}q(A)\Lambda(dA,ds),
    \qquad j\in\bZ^d,
\end{equation}
has a finite fourth moment and variance $\tau^2>0$. The variables $(\varepsilon_j)_{j\in\bZ^d}$ are then IID and centered. A centered compound-Poisson L\'{e}vy basis with bounded jumps provides a concrete specification. This homogeneous input is of the same type as that used by \cite{curato2022}; here non-stationarity is introduced through the $i$- and $n$-dependent kernel and ambit geometry below.

For $v=(v_0,\ldots,v_{d_s})\in[0,1]^d$, fix $\delta\in(0,1)$ and $c_0,c_1,\lambda_0>0$, and define
\[
    a(v)=1+\delta v_0,
    \qquad
    c(v)=c_0+c_1v_1.
\]
Unlike the common lag kernels and translation-invariant ambit sets in \cite{curato2022}, both the loading $a(v)$ and the cone aperture $c(v)$ vary with the scaled location $v=i/n$.
The causal lattice ambit cone and its kernel are
\begin{align*}
    \cA(v)
    &=\left\{u=(r,z)\in\bN_0\times\bZ^{d_s}:
      \|z\|_\infty\leq c(v)r\right\},\\
    \kappa_v(u)
    &=a(v)\exp\{-\lambda_0\|u\|_\infty\}
      \mathbf 1\{u\in\cA(v)\}.
\end{align*}
Define the following triangular array of random fields:
\begin{equation}\label{example_non-stationary_ambit_model}
\begin{aligned}
    X_{i,n} = \sum_{u\in\cA(i/n)}a(i/n)e^{-\lambda_0\|u\|_\infty}\varepsilon_{i-u},
    \qquad i\in D_n,
\end{aligned}
\end{equation}
The following properties hold.
\begin{enumerate}
    \item[(i)] The series in
    \eqref{example_non-stationary_ambit_model} converges absolutely
    almost surely and in $\bL^4$, and
    \[
        \sup_{n\geq1}\sup_{i\in D_n}\|X_{i,n}\|_4<\infty.
    \]
    Moreover, it admits the causal ambit representation
    \begin{equation}\label{example_non-stationary_ambit_integral_representation}
    \begin{aligned}
        X_{i,n}
        =\int_{\mathcal M}\int_{\bR^d}
        g_{i,n}(A,s)\Lambda(dA,ds),
        \qquad i\in D_n,
    \end{aligned}
    \end{equation}
    with kernel
    \begin{equation}\label{example_non-stationary_ambit_kernel}
        g_{i,n}(A,s)
        =q(A)\sum_{u\in\cA(i/n)}
        a(i/n)e^{-\lambda_0\|u\|_\infty}\mathbf 1\{s\in Q_{i-u}\}.
    \end{equation}
    The resulting triangular array of random fields is centered and non-stationary
    within each row, and its law changes with $n$. In particular,
    \begin{equation}\label{example_non-stationary_ambit_variance}
        \var(X_{i,n})
        =\tau^2a(i/n)^2
        \sum_{u\in\cA(i/n)}
        e^{-2\lambda_0\|u\|_\infty}.
    \end{equation}

    \item[(ii)] Proposition \ref{proposition_BS} applies with $l=0$, $p=4$, and $B_{i,n}(h)=(1+\delta)e^{-\lambda_0h}.$
    Consequently, the array is $\eta$-weakly dependent in the sense of
    Definition \ref{def_of_weak_dep}, with
    \begin{equation}\label{example_non-stationary_ambit_eta}
        \bar\eta(r)
        =\bigO\left(
        r^{d-1}e^{-\lambda_0r/2}
        \right).
    \end{equation}

\end{enumerate}

The verification of the above claims is given in Section \ref{proof_subsection_non-stationary_ambit_mma}.

\section{Asymptotic theory for weakly dependent triangular arrays of random fields}\label{section_limit_theorem}

In this section, we investigate the asymptotic behaviour of a weakly dependent triangular array of random fields $\{X_{i,n}: i \in D_n, n \geq 1\}$ on $D \subset \bR^d\ (d \geq 1)$, which satisfies Assumption \ref{assumption_minimum_distance}. The sequence $(D_n)_{n\geq1}$ consists of finite sampling regions in $D$, and $\lim_{n\to\infty}|D_n|_c = \infty$ represents their expansion as $n \to \infty$.

\begin{assumption}\label{assumption_unif_l1+_bound}
$X_{i,n}$ has finite $l$-th moment for some $l > 1$ uniformly over $D_n$ and $n \geq 1$, i.e.
\begin{equation}\label{unif_l1+_bound}
    \sup_n\sup_{i \in D_n}\expct\left|X_{i,n}\right|^l < \infty.
\end{equation}
\end{assumption}
\begin{remark}
By H\"{o}lder's inequality and Markov's inequality, \eqref{unif_l1+_bound} implies the $\bL^p$ uniform integrability for any $0 < p < l$. i.e.
\begin{equation}\label{unif_l1+_int}
    \lim_{k\to\infty}\sup_n\sup_{i\in D_n}\expct\left[\left|X_{i,n}\right|^{p}\mathbf 1\{\left|X_{i,n}\right| \geq k\}\right] = 0.
\end{equation} See page 216 in \cite{billingsley2008} for the definition.
\end{remark}

\begin{assumption}\label{assumption_epsilon_lln}
The dependence coefficient of $\{X_{i,n}: i \in D_n, n \geq 1\}$ satisfies $\bar\eta_{1,1}(r) = \bigO(r^{-\alpha})$ with $\alpha > d$.
\end{assumption}

Now we are ready to present our LLN as follows. 

\begin{theorem}\label{theorem_lln}
    Let $\{X_{i,n}\in\bR: i \in D_n, n \geq 1\}$ be a triangular array of random fields on $D \subset \bR^d (d \geq 1)$, where $(D_n)_{n\geq 1}$ is a sequence of finite sampling regions in $D$ with $\lim_{n\to\infty}|D_n|_c = \infty$. If Assumption \ref{assumption_minimum_distance}, Assumption \ref{assumption_unif_l1+_bound} and Assumption \ref{assumption_epsilon_lln} are satisfied, then as $n \to \infty$, \[\frac{1}{|D_n|_c}\sum_{i\in D_n}(X_{i,n} - \expct X_{i,n}) \to 0\] in $\bL^1$.
\end{theorem}

Let $S_n = \sum_{i \in D_n}X_{i,n}$ and $\sigma^2_n = \var(S_n).$ We need the following assumptions to state our CLT.

\begin{assumption}\label{assumption_unif_l2+_bound}
$X_{i,n}$ has finite $m$-th moment for some $m > 1$ uniformly over $D_n$ and $n \geq 1$, i.e.
\begin{equation}\label{unif_l2+_bound}
    \sup_n\sup_{i \in D_n}\expct\left|X_{i,n}\right|^m < \infty.
\end{equation}
\end{assumption}

\begin{assumption}\label{assumption_eta_clt}
With the same $m > 2$ in Assumption \ref{assumption_unif_l2+_bound}, the $\eta$-coefficient of $\{X_{i,n}: i \in D_n, n \geq 1\}$ satisfies:

\begin{enumerate}
    \item [(a).] For all $u + v \leq 4$, $\bar\eta_{u,v}(r) = \bigO(r^{-\alpha})$ with $\alpha > \frac{m-1}{m-2}d$;
    \item [(b).] $\bar\eta_{\infty,1}(r) := \sup_u\bar\eta_{u,1}(r) = \bigO(r^{-\beta})$ with $\beta > 2d$.
\end{enumerate}

\end{assumption}

\begin{assumption}\label{assumption_variance}
    $\liminf_{n\to\infty}(|D_n|_c)^{-1}\sigma_n^2 > 0.$
\end{assumption}

Assumption \ref{assumption_variance} is a standard condition in the asymptotic theory literature, as maintained in \cite{bolthausen1982}, \cite{jenish2009}, and \cite{jenish2012}. It is required to prove that $\sigma^2_n$ is asymptotically proportional to $|D_n|_c$ as $n\to\infty$, which ensures that no single summand dominates the sum. Our CLT is given as a theorem below.
\begin{theorem}\label{theorem_clt}
    Let $\{X_{i,n}\in\bR: i \in D_n, n \geq 1\}$ be a triangular array of zero-mean random fields on $D \subset \bR^d (d \geq 1)$, where $(D_n)_{n\in\bN}$ is a sequence of finite sampling regions in $D$ with $\lim_{n\to\infty}|D_n|_c = \infty$. If Assumption \ref{assumption_minimum_distance}, Assumption \ref{assumption_unif_l2+_bound}, Assumption \ref{assumption_eta_clt} and Assumption \ref{assumption_variance} are satisfied, then as $n \to \infty$, \[\sigma_n^{-1}S_n \convd \mathcal N(0,1).\]
\end{theorem}

\begin{remark}
In Section \ref{proof_dynamic_network_estimation} we will show that the example in Section \ref{subsection_dynamic_network_bernoulli_array} satisfies the conditions of Theorems \ref{theorem_lln} and \ref{theorem_clt}. A similar argument verifies that the examples in Sections \ref{subsection_non-stationary_bernoulli_array} and \ref{subsection_non-stationary_ambit_mma} also satisfy these conditions, so we omit the proof for simplicity. The examples in Sections \ref{subsection_non-stationary_bernoulli_array} and \ref{subsection_dynamic_network_bernoulli_array} show that the two theorems cover non-mixing triangular arrays of random fields outside the $\alpha$- and $\phi$-mixing framework of \cite{jenish2009}, while the example in Section \ref{subsection_non-stationary_ambit_mma} extends the stationary random-field results of \cite{curato2022} to non-stationary triangular arrays of random fields.
\end{remark}

Theorem \ref{theorem_clt} is stated for triangular arrays of scalar-valued random fields, whereas vector-valued scores and estimating equations are central to statistical inference. Facilitated by the transformation-invariance of $\eta$-weak dependence in Proposition \ref{proposition_transformation}, Theorem \ref{theorem_clt} extends to triangular arrays of vector-valued random fields through a standard Cram\'{e}r-Wold device.
\begin{corollary}\label{corollary_clt}
Let $\{X_{i,n}\in\bR^k: i \in D_n, n \geq 1\}$ be a triangular array of zero-mean vector-valued random fields. By regarding $|\cdot|$ in Assumption \ref{assumption_unif_l2+_bound} as Euclidean norm, and replacing Assumption \ref{assumption_variance} by \[\liminf_{n\to\infty}(|D_n|_c)^{-1}\lambda_{\min}(\Sigma_n) > 0\] where $\lambda_{\min}(\Sigma_n)$ is the smallest eigenvalue of $\Sigma_n := \var(S_n)$, then as $n\to\infty$:\[\Sigma_n^{-1/2}S_n \convd \mathcal N(0, I_k).\]
\end{corollary}

\section{An application in parameter estimation of the model with dynamic network}\label{section_dynamic_network_estimation}

In this section, we utilize the the heredity results in Section \ref{section_heredity} and the asymptotic theory in Section \ref{section_limit_theorem}, and show how preceding results yield a complete parameter estimation, using the spatio-temporal model with dynamic network in Section \ref{subsection_dynamic_network_bernoulli_array} as an example. Throughout this section, the network threshold $\gamma\in\Gamma$ is fixed and known. Let the observations be generated with true parameters
\[
    \xi_0=(a_0,\beta_0')'
    \in\Xi=[\underline a,\bar a]\times\cB,
\]
and write $X_{i,t,n}=X_{i,t,n}(a_0,\beta_0,\gamma)$.

Using the backward-path notation in \eqref{example_dynamic_network_ancestor},
for $(i,t)\in D_n$ and $\xi=(a,\beta')'\in\Xi$, let
\begin{equation}\label{estimation_dynamic_network_mean}
\begin{aligned}
    \mu_{i,t,n}(\xi)
    &=aX_{A_1^\gamma(i,t),t-1,n}+p_{i,t,n}(\beta),\\
    e_{i,t,n}(\xi)&=X_{i,t,n}-\mu_{i,t,n}(\xi).
\end{aligned}
\end{equation}
Innovation independence gives
$\expct(X_{i,t,n}\mid\cF_{t-1,n})=\mu_{i,t,n}(\xi_0)$, where
$\cF_{t-1,n}$ is the $\sigma$-algebra generated by historical innovations up to $t-1$.
The nonlinear least-squares criterion and estimator are
\begin{equation}\label{estimation_dynamic_network_nls}
    Q_n(\xi)=\frac{1}{2|D_n|_c}\sum_{(i,t)\in D_n}e_{i,t,n}(\xi)^2,
    \qquad
    \widehat\xi_n\in\argmin_{\xi\in\Xi}Q_n(\xi).
\end{equation}
Thus \eqref{estimation_dynamic_network_nls} is a conditional-mean estimator.

Let $p^0_{i,t,n}=p_{i,t,n}(\beta_0)$ and define
\begin{equation}\label{estimation_dynamic_network_information}
\begin{aligned}
    J_n&=\frac{1}{|D_n|_c}\sum_{(i,t)\in D_n}
       \expct\left\{
       \nabla_\xi\mu_{i,t,n}(\xi_0)
       \nabla_\xi\mu_{i,t,n}(\xi_0)'\right\},\\
    \Omega_n&=\frac{1}{|D_n|_c}\sum_{(i,t)\in D_n}
       p^0_{i,t,n}(1-p^0_{i,t,n})
       \expct\left\{
       \nabla_\xi\mu_{i,t,n}(\xi_0)
       \nabla_\xi\mu_{i,t,n}(\xi_0)'\right\}.
\end{aligned}
\end{equation}

\begin{assumption}\label{assumption_dynamic_network_estimation}
The following conditions hold.
\begin{enumerate}
    \item[(a)] The true value $\xi_0$ is in the interior of the compact set $\Xi$, and $\sup_{i,t,n}\|z_{i,t,n}\|<\infty$.
    \item[(b)] For every $\delta>0$,
    \begin{equation}\label{estimation_dynamic_network_identification}
        \liminf_{n\to\infty}
        \inf_{\substack{\xi\in\Xi\\\|\xi-\xi_0\|\geq\delta}}
        \frac{1}{|D_n|_c}\sum_{(i,t)\in D_n}
        \expct\bigl[\mu_{i,t,n}(\xi)-\mu_{i,t,n}(\xi_0)\bigr]^2>0.
    \end{equation}
    \item[(c)] $J_n\to J$ and $\Omega_n\to\Omega$, where $J$ and $\Omega$ are positive definite matrices.
\end{enumerate}
\end{assumption}

The first step is to transfer the weak dependence of $X_{i,t,n}$ established in Section \ref{proof_example_dynamic_network_bernoulli_array} to the criterion contributions and their derivatives. Each of these statistics is a function of the local vector
\begin{equation}\label{estimation_dynamic_network_local_vector}
    V_{i,t,n}=\bigl(X_{i,t,n},X_{A_1^\gamma(i,t),t-1,n},z_{i,t,n}'\bigr)'.
\end{equation}

\begin{proposition}\label{proposition_dynamic_network_score_wd}
Consider the following scalar triangular arrays of random fields, where
$c\in\bR^{1+d_\beta}$ with $\|c\|\leq1$ and
$k,\ell\in\{1,\ldots,1+d_\beta\}$:
\[
\begin{gathered}
    \frac12e_{i,t,n}(\xi)^2,
    \qquad
    c'\nabla_\xi\!\left\{\frac12e_{i,t,n}(\xi)^2\right\},
    \\
    \left[\nabla^2_{\xi\xi'}\!\left\{\frac12e_{i,t,n}(\xi)^2\right\}\right]_{k\ell},
    \qquad
    \left[
    \nabla_\xi\!\left\{\frac12e_{i,t,n}(\xi)^2\right\}
    \nabla_\xi\!\left\{\frac12e_{i,t,n}(\xi)^2\right\}'
    \right]_{k\ell}.
\end{gathered}
\]
For every listed triangular array of random fields, there is a constant $C<\infty$, independent of
$n$, $\xi$, $c$, $k$, and $\ell$, such that its weak-dependence coefficient
satisfies
\begin{equation}\label{estimation_dynamic_network_statistic_eta}
    \sup_{\xi\in\Xi}\bar\eta^\xi(r)
    \leq C \bar a^{\lceil(r-2)/2\rceil},
    \qquad r\geq3.
\end{equation}
\end{proposition}

Proposition \ref{proposition_dynamic_network_score_wd} follows by first bounding the coefficient of the local-vector triangular array of random fields in \eqref{estimation_dynamic_network_local_vector} and then applying Proposition \ref{proposition_transformation} in its globally Lipschitz case. It makes Theorems \ref{theorem_lln} and \ref{theorem_clt} directly available for the objective, score, Hessian, and covariance estimator.

\begin{theorem}\label{theorem_dynamic_network_estimator}
Suppose that Assumption \ref{assumption_dynamic_network_estimation} holds. Then, as $n\to\infty$,
\begin{enumerate}
    \item[(i)] (Consistency) $\widehat\xi_n\convp\xi_0$;
    \item[(ii)] (Asymptotic normality)
    \begin{equation}\label{estimation_dynamic_network_asymptotic_normality}
        \sqrt{|D_n|_c}(\widehat\xi_n-\xi_0)
        \convd \mathcal N\bigl(0,J^{-1}\Omega J^{-1}\bigr);
    \end{equation}
    \item[(iii)] with
    \begin{equation}\label{estimation_dynamic_network_sandwich}
    \begin{aligned}
        \widehat J_n&=\frac{1}{|D_n|_c}\sum_{(i,t)\in D_n}
            \left.
            \nabla^2_{\xi\xi'}\!\left\{
            \frac12e_{i,t,n}(\xi)^2\right\}
            \right|_{\xi=\widehat\xi_n},\\
        \widehat\Omega_n&=\frac{1}{|D_n|_c}\sum_{(i,t)\in D_n}
            \left[
            \left.\nabla_\xi\!\left\{
            \frac12e_{i,t,n}(\xi)^2\right\}
            \right|_{\xi=\widehat\xi_n}\right]
            \left[
            \left.\nabla_\xi\!\left\{
            \frac12e_{i,t,n}(\xi)^2\right\}
            \right|_{\xi=\widehat\xi_n}\right]',
    \end{aligned}
    \end{equation}
    we have $\widehat J_n\convp J$, $\widehat\Omega_n\convp\Omega$, and
    $\widehat J_n^{-1}\widehat\Omega_n\widehat J_n^{-1}\convp
    J^{-1}\Omega J^{-1}$.
\end{enumerate}
\end{theorem}

The proofs of Proposition \ref{proposition_dynamic_network_score_wd} and Theorem \ref{theorem_dynamic_network_estimator} are given in Section \ref{proof_dynamic_network_estimation}.

\begin{remark}[Estimating the network threshold]
The fixed-$\gamma$ formulation gives regular inference for the smooth parameter $\xi=(a,\beta')'$. If $\gamma$ is unknown, one may profile \eqref{estimation_dynamic_network_nls} over the ordered values of $q_t$. The resulting discontinuous argmin connects to the threshold and change-point methods of \cite{chan1993} and \cite{hansen2000}; the likelihood-ratio inversion of \cite{hansen2000} provides a template for confidence sets. With endogenous regressors or transition variables, the first-differenced GMM approach of \cite{seoshin2016} gives an alternative. Joint spatio-temporal inference for a dynamically switching network is left for future research.
\end{remark}

The argument extends beyond the least-squares criterion in \eqref{estimation_dynamic_network_nls}. For smooth M-, QML-, and GMM contributions that are uniformly Lipschitz functions of network lags with bounded propagation radius and overlap multiplicity, Proposition \ref{proposition_transformation} transfers weak dependence to the score and Hessian, while Theorems \ref{theorem_lln} and \ref{theorem_clt} provide the LLN and CLT. This gives a common probabilistic route to local-network vector autoregressions \citep{zhu2017}, nonlinear network autoregressions \citep{armillotta2023}, and time series on dynamic networks \citep{krampe2019}. Proposition \ref{proposition_BS} adds infinite-memory propagation under its sensitivity-decay condition. The QML and GMM analyses of dynamic spatial panels with time-varying or endogenous weights \citep{leeyu2012,quleeyu2017,kuersteinerprucha2020} identify further settings in which a uniform weak-dependence representation leads to parameter inference.

\section{Summary}\label{section_conclusion}

This paper develops an asymptotic theory for non-stationary triangular arrays of random fields under row-uniform $\eta$-weak dependence. The coefficient formulation accommodates multi-dimensional lattices with expanding finite sampling regions while allowing both the distributions and dependence structures to vary across locations and rows. The heredity results demonstrate that this dependence structure is preserved under locally Lipschitz transformations and infinite-memory Bernoulli shifts whose sensitivity coefficients may vary across locations and rows. Explicit bounds are derived that translate the primitive dependence of innovations and functional sensitivities into polynomial or exponential decay rates for the observed triangular array of random fields.

The model constructions illustrate the breadth and flexibility of the framework. The two autoregressive triangular arrays of random fields driven by Bernoulli innovations, including the specification with regime-dependent network matchings, possess exponentially decaying $\eta$-coefficients despite the fact that their $\alpha$- and $\phi$-mixing coefficients do not decay. The triangular array of random fields admitting non-stationary ambit and mixed moving-average representations further demonstrates that location-dependent kernels and causal ambit geometries can be accommodated within the same theoretical framework while providing directly verifiable moment, dependence, and variance conditions. Collectively, these examples link the abstract coefficient assumptions to concrete mechanisms, including discrete innovations, evolving propagation structures, and heterogeneous spatio-temporal dynamics.

Under uniform moment, coefficient-decay, and variance conditions, the paper establishes an $\bL^1$ law of large numbers and a central limit theorem for the proposed framework. In the dynamic-network setting, these probabilistic results form the foundation of a nonlinear least-squares estimation procedure. Specifically, the local objective function, score vector, Hessian matrix, and covariance components inherit the underlying uniform $\eta$-weak dependence structure. As a consequence, the resulting estimators are shown to be consistent and asymptotically normal, and a consistent sandwich covariance estimator is obtained for statistical inference on the model parameters.

\newpage

\appendix

\section{Proofs of the lemmas, propositions, and theorems}\label{appendix_proofs}

\begin{lemma}\label{lma_cov_ineq}
Let $\{X_{i,n}: i\in D_n, n\geq 1\}$ be an $\eta$-weakly dependent triangular array of $\bR$-valued random fields. If $\sup_n\sup_{i \in D_n}\norm{X_{i,n}}_p < \infty$ for some $p > 2$, then for any $i, j \in D_n$:
\begin{equation}\label{cov_ineq}
    |\cov(X_{i,n}, X_{j,n})| \leq C \|X\|_p^{\frac{p}{p-1}}[\eta_{n,1,1}(\rho(i,j))]^{\frac{p-2}{p-1}},
\end{equation} where $\|X\|_p := \sup_n\sup_{i \in D_n}\norm{X_{i,n}}_p$.
\end{lemma}

\begin{proof}

Let $X_{i,n}(k) = -k\vee X_{i,n}\wedge k$ be a truncation of $X_{i,n}$ at level $k > 0$, where $\vee$ and $\wedge$ mean \textit{maximum} and \textit{minimum} respectively. Then for any $i\in D_n$ and $a \in (0, p)$:
\begin{align*}
    \expct|X_{i,n} - X_{i,n}(k)|^a \leq &\expct\left[|X_{i,n}|^a\mathbf 1\{|X_{i,n}| \geq k\}\right]\\
    \leq &\left(\expct|X_{i,n}|^p\right)^{a/p}[\prob(|X_{i,n}| \geq k)]^{1-(a/p)}\\
    \leq &\|X\|_p^a\left[\frac{\|X\|_p^p}{k^p}\right]^{1-(a/p)} = \|X\|_p^p k^{a-p},
\end{align*} where the second line and the third line come from H\"{o}lder's inequality and Markov inequality respectively. Hence $\sup_n\sup_{i\in D_n}\|X_{i,n} - X_{i,n}(k)\|_a \leq \|X\|_p^{p/a}k^{1-(p/a)}$.

For any $i, j\in D_n$ we have
\begin{align*}
    |\cov(X_{i,n}, X_{j,n})| \leq &|\cov(X_{i,n}(k), X_{j,n}(k))|\\
    &+ |\cov(X_{i,n} - X_{i,n}(k), X_{j,n}(k))|\\
    &+ |\cov(X_{i,n}, X_{j,n} - X_{j,n}(k))|.
\end{align*} For the last term on the right-hand-side (RHS), we could find $a \in (1, p)$ such that $1/a + 1/p = 1$ and therefore
\begin{align*}
    &|\cov(X_{i,n}, X_{j,n} - X_{j,n}(k))|\\
    \leq &|\expct[X_{i,n}(X_{j,n} - X_{j,n}(k))]| + |\expct(X_{i,n})||\expct(X_{j,n} - X_{j,n}(k))|\\
    \leq &\|X_{i,n}\|_p\|X_{j,n} - X_{j,n}(k)\|_a + \|X_{i,n}\|_1\|X_{j,n} - X_{j,n}(k)\|_1\\
    \leq &2\|X_{i,n}\|_p\|X_{j,n} - X_{j,n}(k)\|_a\\
    \leq &2 \|X\|_p^p k^{2-p}.
\end{align*} Same bound could also be derived for the second term on the RHS. As for the first term, note that $X(k)$ is a function of $X$ with bound $k$ and Lipschitz constant 1. Then by \eqref{cov_ineq_bounded} we have $|\cov(X_{i,n}(k), X_{j,n}(k))| \leq 2k\eta_{n,1,1}(\rho(i,j))$, then \[|\cov(X_{i,n}, X_{j,n})| \leq 6\|X\|_p^{\frac{p}{p-1}}[\eta_{n,1,1}(\rho(i,j))]^{\frac{p-2}{p-1}}\] by choosing $k = \left[\frac{\|X\|_p^p}{\eta_{n,1,1}(\rho(i,j))}\right]^{\frac{1}{p-1}}$.

\end{proof}

\begin{lemma}\label{lma_brockwell}
    (Proposition 6.3.9 in \cite{brockwell2009}) Let $(Z_n)_{n \geq 1}$ and $(V_{n,k})_{n,k \in \bN}$ be sequences of random vectors. $Z_n \convd V$ if the following statements are true:
    \begin{enumerate}
        \item For each $k \in \bN$, there exists $V_k$ such that $V_{n,k} \convd V_k$ as $n \to \infty$;
        \item $V_k \convd V$ as $k \to \infty$;
        \item $\lim_{k \to \infty}\limsup_{n \to \infty}\prob(|Z_n - V_{n,k}| > \delta) = 0$ for any $\delta > 0$.
    \end{enumerate}
\end{lemma}

\begin{lemma}\label{lma_stein}
    (Lemma 2 in \cite{bolthausen1982}) Let $(\nu_n)_{n\in\bN}$ be a sequence of probability measures over $\bR$ with
    \begin{enumerate}
        \item $\sup_n\int x^2\nu_n(dx) < \infty$,
        \item $\lim_{n\to\infty}\int(\img\lambda - x)e^{\img\lambda x}\nu_n(dx) = 0$ for all $\lambda\in\bR$.
    \end{enumerate}
    Then $\nu_n \convd \mathcal N(0,1)$ as $n\to\infty$.
\end{lemma}

\begin{lemma}\label{lma_jenish2009}
    (Lemma A.1 in \cite{jenish2009}) For any $i \in D\subset\bR^d$ and $h \geq 1$, let \[N_i(h) := \left|\left\{j \in D: h \leq \rho(i,j) < h+1\right\}\right|_c\] be the number of all elements of $D$ located at any distance in $[h, h+1)$ from $i$. Then $\sup_i N_i(h) \leq C h^{d-1}$.
\end{lemma}

\begin{lemma}\label{lma_xu2022}
    (Lemma A.4 in \cite{xu2022}) For any $\alpha > 0$ and $s \geq 2$, \[\sum_{h=\lfloor s\rfloor}^{\infty}h^{-\alpha-1} < \frac{2^{\alpha+1}}{\alpha}s^{-\alpha}.\]
\end{lemma}

\subsection{Proof of results in Section \ref{section_weak_dependence}}

\subsubsection{Verification of the example in Section \ref{subsection_non-stationary_bernoulli_array}}\label{proof_example_non-stationary_bernoulli_array}

We first verify part (i). The series in \eqref{example_bernoulli_array} converges absolutely, is bounded between zero and one, and is causal because it depends only on current and past innovations. Moreover,
\[
\begin{aligned}
\frac12X_{i,t-1,n}+\frac12\varepsilon_{i,t,n}
&=\frac12\sum_{k=0}^{\infty}2^{-(k+1)}
  \varepsilon_{i,t-1-k,n}
  +\frac12\varepsilon_{i,t,n}\\
&=\sum_{k=0}^{\infty}2^{-(k+1)}
  \varepsilon_{i,t-k,n}
=X_{i,t,n},
\end{aligned}
\]
so the series solves the autoregression. If $X_{i,t,n}$ and $\widetilde X_{i,t,n}$ are two bounded solutions, let $M<\infty$ be a common bound. Iterating their difference backward $m$ times gives
\[
    |X_{i,t,n}-\widetilde X_{i,t,n}|
    =2^{-m}|X_{i,t-m,n}-\widetilde X_{i,t-m,n}|
    \leq2M\,2^{-m}.
\]
Letting $m\to\infty$ proves uniqueness. 

We next verify non-stationarity. The series in \eqref{example_bernoulli_array} is the binary expansion
\[
X_{i,t,n}=0.\varepsilon_{i,t,n}\varepsilon_{i,t-1,n}\varepsilon_{i,t-2,n}\cdots\quad\text{(base 2)}.
\]
Because the success probabilities are bounded between $1/4$ and $3/4$, the probability that its digits are eventually all zero or eventually all one is zero. The binary expansion is therefore unique almost surely, and
\begin{equation}\label{proof_example_marginal_prob}
    \prob\left(X_{i,t,n}\geq\frac12\right)=p_{i,t,n}.
\end{equation}
According to \eqref{example_bernoulli_prob}, we have $p_{n,4n,n}=3/4$ and $p_{2n,4n,n}=1/2$, which means that $X_{n,4n,n}$ and $X_{2n,4n,n}$ have different distribution laws, proving within-row non-stationarity.

Apparently, $p_{1,1,n}$ changes with $n$. Moreover, since $p_{2n,t,n}=1/2$ for every $t$, we have $\var(X_{2n,4n,n})=1/12$. By contrast, $p_{n,4n,n}=3/4$ and all the other innovation variances are at most $1/4$, so $\var(X_{n,4n,n})<1/12$. Hence deterministic centering does not restore stationarity. (i) is proved.

We next prove the uniform $\eta$ bound in (ii). Fix $r\geq1$, set $q=\lceil r\rceil-1$, and define
\[
X^{(q)}_{i,t,n}=\sum_{k=0}^{q}2^{-(k+1)}\varepsilon_{i,t-k,n}.
\]
The omitted tail is bounded deterministically by
\begin{equation}\label{proof_example_truncation_bound}
    \left|X_{i,t,n}-X^{(q)}_{i,t,n}\right|
    \leq\sum_{k=q+1}^{\infty}2^{-(k+1)}
    =2^{-(q+1)}=2^{-\lceil r\rceil}.
\end{equation}
Let $U_n,V_n\subset D_n$ have cardinalities $u,v$ and satisfy $\rho(U_n,V_n)\geq r$. If $(i,t)\in U_n$ and $(j,s)\in V_n$ satisfy $i=j$, then $|t-s|\geq\lceil r\rceil=q+1$, so the innovation windows entering $X^{(q)}_{i,t,n}$ and $X^{(q)}_{j,s,n}$ are disjoint. If $i\neq j$, the corresponding innovations are independent by construction. For any $A\subset D_n$, write
\[
    \fX_A^{(q)}
    =\bigl(X_{i,t,n}^{(q)}\bigr)_{(i,t)\in A}.
\] It follows that the two truncated vectors $\fX^{(q)}_{U_n}$ and $\fX^{(q)}_{V_n}$ are independent.

For arbitrary $f\in\cF_u$ and $g\in\cG_v$, independence of the truncated vectors and the elementary inequality $|\cov(A,B)|\leq2\|B\|_{\infty}\expct|A|$ give
\begin{align*}
&\left|\cov\left(f(\fX_{U_n}),g(\fX_{V_n})\right)\right|\\
\leq{}&2\|g\|_{\infty}\expct\left|f(\fX_{U_n})-f(\fX^{(q)}_{U_n})\right|
+2\|f\|_{\infty}\expct\left|g(\fX_{V_n})-g(\fX^{(q)}_{V_n})\right|\\
\leq{}&2^{1-\lceil r\rceil}
\left[u\|g\|_{\infty}\lip(f)+v\|f\|_{\infty}\lip(g)\right]\\
={}&2^{1-\lceil r\rceil}\Psi(f,g),
\end{align*}
where \eqref{proof_example_truncation_bound} is used in the second inequality. Dividing by $\Psi(f,g)$ and taking the suprema over $f,g,U_n,V_n,n,u,v$ proves \eqref{example_eta_bound}.

It remains to verify the failure of mixing. Fix an integer $r\geq1$ and choose $n$ such that $r+1\leq4n$. Both $(1,1)$ and $(1,1+r)$ then belong to $D_n$.

By uniqueness of the binary expansions, the event
\[
A_{n,r}=\{\varepsilon_{1,1,n}=1\}
\]
is measurable, up to a null set, with respect to both $\sigma(X_{1,1,n})$ and $\sigma(X_{1,1+r,n})$: it is the first binary digit of the former variable and the $(r+1)$st binary digit of the latter. Writing $p=p_{1,1,n}\in[1/4,3/4]$ and using the notation of Definition 1 in \cite{jenish2009}, we obtain
\begin{align*}
\alpha_n(\{(1,1)\},\{(1,1+r)\})
&\geq\left|\prob(A_{n,r})-\prob(A_{n,r})^2\right|
=p(1-p)\geq\frac{3}{16},\\
\phi_n(\{(1,1)\},\{(1,1+r)\})
&\geq\left|\prob(A_{n,r}\mid A_{n,r})-\prob(A_{n,r})\right|
=1-p\geq\frac14.
\end{align*}
The two singleton index sets are at distance $r$. Taking the supremum over $n$ yields \eqref{example_mixing_lower_bound}, proving part (iii).
This completes the verification.

\subsubsection{Verification of the example in Section \ref{subsection_dynamic_network_bernoulli_array}}\label{proof_example_dynamic_network_bernoulli_array}

The proof of (i) is similar to the proof of (i) in Section \ref{subsection_non-stationary_bernoulli_array}, therefore omitted here. We next establish the uniform $\eta$ bound. Fix $r\geq1$, let $m=\lceil r/2\rceil$, and define the truncated variable
\[
    X_{i,t,n}^{(m)}
    =\sum_{k=0}^{m-1}a^k
    \varepsilon_{A_k^\gamma(i,t),t-k,n}.
\]
The tail satisfies
\begin{equation}\label{proof_dynamic_network_truncation_bound}
\begin{aligned}
    \left|X_{i,t,n}-X_{i,t,n}^{(m)}\right|
    &\leq\frac{a^m}{1-a}
    \leq\frac{\bar a^m}{1-\bar a}
    =:\delta_m.
\end{aligned}
\end{equation}
For $z=(i,t)$ and $z'=(j,s)$, if $X_{i,t,n}^{(m)}$ and $X_{j,s,n}^{(m)}$ share a common innovation, then there exist $k,l\in\{0,\ldots,m-1\}$ such that
\begin{equation}\label{proof_dynamic_network_common_innovation}
    \bigl(A_k^\gamma(i,t),t-k\bigr)
    =\bigl(A_l^\gamma(j,s),s-l\bigr).
\end{equation}
According to \eqref{example_dynamic_network_matchings}, any two neighboring nodes' indices have difference one, therefore 
\[
    \bigl|A_k^\gamma(i,t)-i\bigr|\leq k,
    \qquad
    \bigl|A_l^\gamma(j,s)-j\bigr|\leq l.
\]
By \eqref{proof_dynamic_network_common_innovation} we have
\[
    |i-j|\leq \bigl|i-A_k^\gamma(i,t)\bigr| + \bigl|A_l^\gamma(j,s)-j\bigr| \leq k+l\leq 2m-2,
\]
and
\[
    |t-s|=|k-l|\leq m-1.
\]
Since $2\lceil r/2\rceil-2<r$, this implies $\rho(z,z')<r$. Consequently, if $U_n,V_n\subset D_n$ satisfy $\rho(U_n,V_n)\geq r$, the two truncated vectors $\fX_{U_n}^{(m)}$ and $\fX_{V_n}^{(m)}$ depend on disjoint sets of innovations and are therefore independent.

For arbitrary bounded Lipschitz functions $f\in\cF_u$ and $g\in\cG_v$, this independence and \eqref{proof_dynamic_network_truncation_bound} yield
\begin{align*}
&\left|\cov\left(f(\fX_{U_n}),g(\fX_{V_n})\right)\right|\\
\leq{}&2\|g\|_\infty
    \expct\left|f(\fX_{U_n})-f(\fX_{U_n}^{(m)})\right|
    +2\|f\|_\infty
    \expct\left|g(\fX_{V_n})-g(\fX_{V_n}^{(m)})\right|\\
\leq{}&2\delta_m
    \left[u\|g\|_\infty\lip(f)
    +v\|f\|_\infty\lip(g)\right]
 =2\delta_m\Psi(f,g).
\end{align*}
Dividing by $\Psi(f,g)$ proves \eqref{example_dynamic_network_eta_bound}.

We now prove the failure of mixing. The restriction $a<1/2$ separates the conditional supports in \eqref{example_dynamic_network_model}: if $\varepsilon_{i,t,n}=0$, then
\[
    0\leq X_{i,t,n}\leq\frac{a}{1-a}<1,
\]
whereas $\varepsilon_{i,t,n}=1$ implies $X_{i,t,n}\geq1$. Hence:
\begin{equation}\label{proof_dynamic_network_decoder}
\begin{split}
    \varepsilon_{i,t,n}
    &=\mathbf 1\{X_{i,t,n}\geq1\},\\
    \mathcal I_a(X_{i,t,n})
    &=X_{\pi_t^\gamma(i),t-1,n},
\end{split}
\end{equation} where $\mathcal I_a(x)=\frac{x-\mathbf 1\{x\geq1\}}{a}$.

Fix an integer $r\geq1$ and choose any $n\geq r$. Starting from $d_0=2n$, define the unique forward network path by
\begin{equation}\label{proof_dynamic_network_descendant_path}
    d_h=(\pi_{n+h}^\gamma)^{-1}(d_{h-1}),
    \qquad 1\leq h\leq r.
\end{equation}
Because each step changes the node index by one, $|d_h-2n|\leq h$. Therefore both $z_0=(2n,n)$ and $z_r=(d_r,n+r)$ belong to $D_n$, and
\[
    \rho(z_0,z_r)=r.
\]
Repeated application of \eqref{proof_dynamic_network_decoder} gives
\begin{equation}\label{proof_dynamic_network_recursive_decoder}
    \mathcal I_a^r\bigl(X_{d_r,n+r,n}\bigr)
    =X_{2n,n,n}.
\end{equation}
Here $\mathcal I_a^r$ denotes the $r$-fold composition of $\mathcal I_a$.
The event
\[
    \{X_{2n,n,n}\geq1\}
    =\{\varepsilon_{2n,n,n}=1\}
\]
belongs to $\sigma(X_{2n,n,n})$. By \eqref{proof_dynamic_network_recursive_decoder}, the almost surely identical event
\[
    \left\{\mathcal I_a^r
    \bigl(X_{d_r,n+r,n}\bigr)\geq1\right\}
\]
belongs to $\sigma(X_{d_r,n+r,n})$. Writing $p=p_{2n,n,n}(\beta)\in[p_0,1-p_0]$, the definitions of the singleton coefficients give
\begin{align*}
\alpha_n^\vartheta(\{z_0\},\{z_r\})
&\geq p-p^2
=p(1-p)\geq p_0(1-p_0),\\
\phi_n^\vartheta(\{z_0\},\{z_r\})
&\geq1-p\geq p_0.
\end{align*}
Taking the supremum over $n$ and then the infimum over $\vartheta$ proves \eqref{example_dynamic_network_mixing_lower_bound}.
This completes the verification.

\subsubsection{Proof of Proposition \ref{proposition_BS}}

Firstly we prove that $X_{i,n} = H_{i,n}((\varepsilon_j)_{j\in D})$ is well-defined in $\bL^1$. For any $s \in \bN_0$, let $X_{i,n}^{(s)} = H_{i,n}((\varepsilon_j\mathbf 1\{\rho(i,j)\leq s\})_{j\in D})$. Then by \eqref{lipschitz_condition_BS} we have
\begin{align*}
    &\left|X_{i,n}^{(s+m)} - X_{i,n}^{(s)}\right|\\
    \leq &\sum_{k=1}^m\left|X_{i,n}^{(s+k)} - X_{i,n}^{(s+k-1)}\right|\\
    = &\sum_{k=1}^m\left|H_{i,n}((\varepsilon_j\mathbf 1\{\rho(i,j)\leq s+k\})_{j\in D}) - H_{i,n}((\varepsilon_j\mathbf 1\{\rho(i,j)\leq s+k-1\})_{j\in D})\right|\\
    \leq &\sum_{k=1}^m B_{i,n}(s+k)(\max_{\rho(i,j)\leq s+k-1}|\varepsilon_j|^l \vee 1)\sum_{\rho(i,j)=s+k}|\varepsilon_j|.
\end{align*} Since $\sup_{i\in D}\|\varepsilon_i\|_p < \infty$ with $p > l + 1$, by H\"{o}lder's inequality and Lemma \ref{lma_jenish2009} we obtain that
\begin{equation}\label{proposition_BS_eq1}
    \norm{X_{i,n}^{(s+m)} - X_{i,n}^{(s)}}_1 \leq C\sum_{k=1}^m (s+k)^{d-1}B_{i,n}(s+k).
\end{equation} Notice that $(s+k)^{d-1}B_{i,n}(s+k) \rightarrow 0$ as $s\to\infty$, according to \eqref{coefficient_condition_BS}. Then if $m$ is fixed, $\norm{X_{i,n}^{(s+m)} - X_{i,n}^{(s)}}_1 \rightarrow 0$ as $s\to\infty$. Therefore $\{X_{i,n}^{(s)}: s\geq 0\}$ is a Cauchy sequence in $\bL^1$, and $X_{i,n} = \lim_{s\to\infty}X_{i,n}^{(s)}$ is well-defined.

Let $U_n, V_n \subset D_n$ be two arbitrary finite subsets of $D_n$ with $|U_n|_c = u$, $|V_n|_c = v$ and $\rho(U_n, V_n) \geq r$. $f \in \cF_u$ and $g \in \cG_v$ are two arbitrary Lipschitz functions with $\|f\|_{\infty} = \|g\|_{\infty} = 1$. For an arbitrary threshold value $T > 0$, define $\varepsilon_i(T) := -T \vee \varepsilon_i \wedge T$, and $X_{i,n}^{(s)}(T) := H_{i,n}((\varepsilon_j(T)\mathbf 1\{\rho(i,j)\leq s\})_{j\in D})$. Notice that
\begin{equation}\label{proposition_BS_eq2}
\begin{aligned}
    &\left|\cov[f((X_{i,n})_{i\in U_n}), g((X_{i,n})_{i\in V_n})]\right|\\
    \leq &\left|\cov[f((X_{i,n})_{i\in U_n}) - f((X_{i,n}^{(s)}(T))_{i\in U_n}), g((X_{i,n})_{i\in V_n})]\right|\\
    &+\left|\cov[f((X_{i,n}^{(s)}(T))_{i\in U_n}), g((X_{i,n})_{i\in V_n}) - g((X_{i,n}^{(s)}(T))_{i\in V_n})]\right|\\
    &+\left|\cov[f((X_{i,n}^{(s)}(T))_{i\in U_n}), g((X_{i,n}^{(s)}(T))_{i\in V_n})]\right|.
\end{aligned}
\end{equation}

We start with the first term in the right-hand-side (RHS) of \eqref{proposition_BS_eq2}:
\begin{align*}
    &\left|\cov[f((X_{i,n})_{i\in U_n}) - f((X_{i,n}^{(s)}(T))_{i\in U_n}), g((X_{i,n})_{i\in V_n})]\right|\\
    \leq &2\lip(f)\sum_{i\in U_n}\expct|X_{i,n} - X_{i,n}^{(s)}(T)|\\
    \leq &2u\lip(f)\left[\sup_{i\in U_n}\expct|X_{i,n} - X_{i,n}^{(s)}| + \sup_{i\in U_n}\expct|X_{i,n}^{(s)} - X_{i,n}^{(s)}(T)|\right].
\end{align*} Same as \eqref{proposition_BS_eq1}, we can prove that \[\expct|X_{i,n} - X_{i,n}^{(s)}| \leq C\sum_{h\geq s} h^{d-1}B_{i,n}(h).\] Notice that $\sum_{\rho(i,j)\leq s} = \bigO(s^d)$ according to Lemma A.1(ii) in \cite{jenish2009}, then by using \eqref{lipschitz_condition_BS} repeatedly we can also prove that (ignoring a constant factor):
\begin{align*}
    &\expct|X_{i,n}^{(s)} - X_{i,n}^{(s)}(T)|\\
    = &\expct\left|H_{i,n}((\varepsilon_j\mathbf 1\{\rho(i,j)\leq s\})_{j\in D}) - H_{i,n}((\varepsilon_j(T)\mathbf 1\{\rho(i,j)\leq s\})_{j\in D})\right|\\
    \leq &\left(\sum_{h=0}^{\infty}B_{i,n}(h)h^{d-1}\right)\expct\left|(\max_{\rho(i,j)\leq s}|\varepsilon_j|^l)\sum_{\rho(i,j)\leq s}|\varepsilon_j|\mathbf 1\{|\varepsilon_j|\geq T\}\right|\\
    \leq &\left(\sum_{h=0}^{\infty}B_{i,n}(h)h^{d-1}\right)s^d\expct\left|\max_{\rho(i,j)\leq s}|\varepsilon_j|^{l+1}\mathbf 1\{|\varepsilon_j|\geq T\}\right|.
\end{align*} Since $\expct|\varepsilon_i|^p < \infty$ and $p > l+1$, by H\"{o}lder's inequality we have \[\expct|X_{i,n}^{(s)} - X_{i,n}^{(s)}(T)| \leq \left(\sum_{h=0}^{\infty}B_{i,n}(h)h^{d-1}\right) C s^dT^{l+1-p}.\] Then we obtain the bound
\begin{equation}\label{proposition_BS_eq3}
\begin{aligned}
    &\left|\cov[f((X_{i,n})_{i\in U_n}) - f((X_{i,n}^{(s)}(T))_{i\in U_n}), g((X_{i,n})_{i\in V_n})]\right|\\
    \leq &C u\lip(f)\left[R_B(s) + C_B s^dT^{l+1-p}\right].
\end{aligned}
\end{equation} The bound of the second term on RHS of \eqref{proposition_BS_eq2} follows analogously:
\begin{equation}\label{proposition_BS_eq4}
\begin{aligned}
    &\left|\cov[f((X_{i,n}^{(s)}(T))_{i\in U_n}), g((X_{i,n})_{i\in V_n}) - g((X_{i,n}^{(s)}(T))_{i\in V_n})]\right|\\
    \leq &C v\lip(g)\left[R_B(s) + C_B s^dT^{l+1-p}\right].
\end{aligned}
\end{equation}

To obtain the bound for the last term on RHS of \eqref{proposition_BS_eq2}, define following functions $F_T: \bR^{s^du}\mapsto\bR$ and $G_T: \bR^{s^dv}\mapsto\bR$ as follows:
\begin{align*}
    &F_T(((\varepsilon_j)_{\rho(i,j)\leq s})_{i\in U_n}) := f((H_{i,n}((\varepsilon_j(T)\mathbf 1\{\rho(i,j)\leq s\})_{j\in D}))_{i\in U_n}) = f((X_{i,n}^{(s)}(T))_{i\in U_n});\\
    &G_T(((\varepsilon_j)_{\rho(i,j)\leq s})_{i\in V_n}) := g((H_{i,n}((\varepsilon_j(T)\mathbf 1\{\rho(i,j)\leq s\})_{j\in D}))_{i\in V_n}) = g((X_{i,n}^{(s)}(T))_{i\in V_n}).
\end{align*} By the $\eta$-weak dependence of $\{\varepsilon_i: i\in D\}$, if $r \geq 2s$ we have
\begin{align*}
    &\left|\cov[f((X_{i,n}^{(s)}(T))_{i\in U_n}), g((X_{i,n}^{(s)}(T))_{i\in V_n})]\right|\\
    \leq &[s^du\lip(F_T) + s^dv\lip(G_T)]\bar\eta^\varepsilon(r-2s).
\end{align*} Notice that for any $\fX = ((X_{ij,n})_{\rho(i,j)\leq s})_{i\in U_n}, \fY = ((Y_{ij,n})_{\rho(i,j)\leq s})_{i\in U_n} \in \bR^{s^du}$, we have:
\begin{align*}
    &\frac{\left|F_T(\fX) - F_T(\fY)\right|}{\sum_{i\in U_n}|(X_{ij,n})_{\rho(i,j)\leq s}-(Y_{ij,n})_{\rho(i,j)\leq s}|}\\
    \leq &\lip(f)\frac{\sum_{i\in U_n}\left|H_{i,n}((X_{ij,n}(T)\mathbf 1\{\rho(i,j)\leq s\})_{j\in D}) - H_{i,n}((Y_{ij,n}(T)\mathbf 1\{\rho(i,j)\leq s\})_{j\in D})\right|}{\sum_{i\in U_n}\sum_{\rho(i,j)\leq s}|X_{ij,n}-Y_{ij,n}|}\\
    \leq &\lip(f)\left(\sum_{h=0}^{\infty}B_{i,n}(h)h^{d-1}\right)T^l\frac{\sum_{i\in U_n}\sum_{\rho(i,j)\leq s}|X_{ij,n}(T)-Y_{ij,n}(T)|}{\sum_{i\in U_n}\sum_{\rho(i,j)\leq s}|X_{ij,n}-Y_{ij,n}|}
\end{align*} by using \eqref{lipschitz_condition_BS} repeatedly. Therefore we can prove that:
\begin{align*}
    &\lip(F_T) \leq C_B T^l\lip(f);\\
    &\lip(G_T) \leq C_B T^l\lip(g).
\end{align*} Then we can bound the last term on RHS of \eqref{proposition_BS_eq2} by
\begin{equation}\label{proposition_BS_eq5}
\begin{aligned}
    &\left|\cov[f((X_{i,n}^{(s)}(T))_{i\in U_n}), g((X_{i,n}^{(s)}(T))_{i\in V_n})]\right|\\
    \leq &[u\lip(f) + v\lip(g)] C_B s^d T^l \bar\eta^\varepsilon(r-2s).
\end{aligned}
\end{equation} for any $r \geq 2s$.

Combining \eqref{proposition_BS_eq3}, \eqref{proposition_BS_eq4}, \eqref{proposition_BS_eq5} we can prove Proposition \ref{proposition_BS} by setting the threshold value $T =\bar\eta^\varepsilon(r-2s)^{-\frac{1}{p-1}}$.

\subsubsection{Verification of \texorpdfstring{\eqref{eq_wd_bound_power}}{} and \texorpdfstring{\eqref{eq_wd_bound_exp}}{}}\label{eq_wd_bound_proof}

Assume that $\bar\eta^\varepsilon(r) = \bigO(r^{-\mu})$ for some $\mu > \frac{p-1}{p-1-l}d$ and $B_{i,n}(h) = \bigO(h^{-b})$ for some $b \geq \frac{p-1-l}{p-1}\mu$. Notice that:
\[\sum_{h=s}^{\infty}h^{-k} \leq \int_{s-1}^{\infty}\frac{1}{x^k}dx = \frac{1}{k-1}(s-1)^{1-k}\] if $k > 1$. Then we have $R_B(s) \leq C s^{-b+d}$. \eqref{eq_wd_bound_power} follows by letting $s = \lfloor r/3\rfloor < \lfloor r/2\rfloor$ in \eqref{eq_eta_output}.

Assume that $d=2$, $\bar\eta^\varepsilon(r) = \bigO(r^{-\mu})$ for some $\mu > 0$ and $B_{i,n}(h) = \bigO(e^{-bh})$ for some $b \geq \frac{p-1-l}{p-1}\mu$. Notice that: 
\[\sum_{h=s}^{\infty} he^{-bh} = \frac{se^{-bs} - (s-1)e^{-b(s+1)}}{(1-e^{-b})^2}.\]
Then $R_B(s) = \bigO(se^{-bs})$, and \eqref{eq_wd_bound_exp} follows by letting $s = \lfloor\log r\rfloor < \lfloor r/2\rfloor$ in \eqref{eq_eta_output}.

\subsubsection{Verification of the claims in Section \ref{subsection_non-stationary_ambit_mma}}\label{proof_subsection_non-stationary_ambit_mma}

The integrals of the L\'{e}vy basis $\Lambda$ over the disjoint cells $(Q_j)_{j\in\bZ^d}$ in \eqref{example_ambit_innovations} are IID. The centering of $\Lambda$ gives $\expct\varepsilon_j=0$, and the moment specification gives $\|\varepsilon_j\|_4<\infty$ and $\var(\varepsilon_j)=\tau^2$. Recall from \eqref{example_non-stationary_ambit_model} that
$$X_{i,n} = \sum_{u\in\cA(i/n)}
      a(i/n)e^{-\lambda_0\|u\|_\infty}\varepsilon_{i-u} = \sum_{u\in\cA(i/n)} \kappa_{i/n}(u)\varepsilon_{i-u}$$

For every $v\in[0,1]^d$,
\begin{equation}\label{proof_ambit_kernel_sum}
    \sum_{u\in\cA(v)}\kappa_v(u)
    \leq(1+\delta)\sum_{u\in\bZ^d}
    e^{-\lambda_0\|u\|_\infty}
    =: C_\kappa < \infty.
\end{equation}
The finiteness of $C_\kappa$ follows from Lemma \ref{lma_jenish2009}.
It follows that
\[
    \expct\left[
      \sum_{u\in\cA(i/n)}\kappa_{i/n}(u)|\varepsilon_{i-u}|
    \right]
    \leq C_\kappa \expct|\varepsilon_0|<\infty.
\]
Hence the series is absolutely convergent almost surely. Minkowski's inequality also gives
\[
    \|X_{i,n}\|_4
    \leq\sum_{u\in\cA(i/n)}
      \kappa_{i/n}(u)\|\varepsilon_{i-u}\|_4
    \leq C_\kappa \|\varepsilon_0\|_4,
\]
uniformly in $i$ and $n$.

To verify \eqref{example_non-stationary_ambit_integral_representation}, let
\[
    \cA_m(v)=\{u\in\cA(v):\|u\|_\infty\leq m\}
\]
and define $g_{i,n}^{(m)}$ by replacing $\cA(i/n)$ with $\cA_m(i/n)$ in \eqref{example_non-stationary_ambit_kernel}. By \eqref{example_ambit_innovations} and linearity of the stochastic integral over the finite collection of cells,
\[
\begin{aligned}
    \int_{\mathcal M}\int_{\bR^d}
    g_{i,n}^{(m)}(A,s)\Lambda(dA,ds)
    &=\sum_{u\in\cA_m(i/n)}
    \kappa_{i/n}(u)\varepsilon_{i-u}.
\end{aligned}
\]
Moreover, Minkowski's inequality and \eqref{proof_ambit_kernel_sum} give
\[
\begin{aligned}
    \left\|
    X_{i,n}-\sum_{u\in\cA_m(i/n)}
    \kappa_{i/n}(u)\varepsilon_{i-u}
    \right\|_4
    &\leq(1+\delta)\|\varepsilon_0\|_4
    \sum_{\substack{u\in\bZ^d\\\|u\|_\infty>m}}
    e^{-\lambda_0\|u\|_\infty}
    \longrightarrow0.
\end{aligned}
\]
The stochastic integral with kernel $g_{i,n}$ is the $\bL^4$ limit of the simple-kernel integrals with kernels $g_{i,n}^{(m)}$. Taking limits in the finite-sum identity therefore proves \eqref{example_non-stationary_ambit_integral_representation}.

Independence and centering of the innovations imply
\[
    \var(X_{i,n})
    =\tau^2\sum_{u\in\cA(i/n)}\kappa_{i/n}(u)^2,
\]
which is \eqref{example_non-stationary_ambit_variance}. Since $a(i/n)$ and $\cA(i/n)$ vary with the scaled coordinates, the array is non-stationary within each row and its law changes with $n$.

We next verify the conditions of Proposition \ref{proposition_BS}. For $x=(x_j)_{j\in\bZ^d}$, define
\[
    H_{i,n}(x)
    =\sum_{u\in\cA(i/n)}\kappa_{i/n}(u)x_{i-u}
\]
on the set where the series is absolutely convergent, and set $H_{i,n}(x)=0$ elsewhere. Every distance shell in $\bZ^d$ is finite. Thus two inputs differing only on the shell $\{j:\rho(i,j)=h\}$ either both belong to the absolute-convergence set or both lie outside it. On the former set,
\begin{align*}
    |H_{i,n}(x)-H_{i,n}(y)|
    &\leq\sum_{\rho(i,j)=h}
      \kappa_{i/n}(i-j)|x_j-y_j|\\
    &\leq(1+\delta)e^{-\lambda_0h}
      \sum_{\rho(i,j)=h}|x_j-y_j|,
\end{align*}
and on the latter set the left-hand side is zero. Therefore \eqref{lipschitz_condition_BS} holds with $l=0$ and
\[
    B_{i,n}(h)=(1+\delta)e^{-\lambda_0h}.
\]
In particular,
\[
    \sup_{n\geq1}\sup_{i\in D_n}
    \sum_{h=0}^{\infty}B_{i,n}(h)h^{d-1}<\infty,
\]
and, for $s\geq1$,
\begin{equation}\label{proof_ambit_sensitivity_tail}
    R_B(s)
    \leq C\sum_{h\geq s}h^{d-1}e^{-\lambda_0h}
    \leq C s^{d-1}e^{-\lambda_0s}.
\end{equation}

The innovations are independent, so $\bar\eta^\varepsilon(r)=0$ for every $r\geq1$. For $r\geq3$, take
\[
    s_r=\left\lfloor\frac{r-1}{2}\right\rfloor.
\]
Then $r-2s_r\geq1$ and $s_r$ is proportional to $r$. Proposition \ref{proposition_BS} and \eqref{proof_ambit_sensitivity_tail} give
\[
    \bar\eta(r)
    \leq C R_B(s_r)
    \leq C r^{d-1}e^{-\lambda_0r/2}.
\]
This proves \eqref{example_non-stationary_ambit_eta}.
This completes the verification.

\subsection{Proof of Theorem \ref{theorem_lln}}

For $k>0$, decompose $X_{i,n}$ as
\begin{equation}\label{truncation}
\begin{array}{ll}
     &X_{i,n}(k) = -k \vee X_{i,n} \wedge k, \\
     &\widetilde{X}_{i,n}(k) = X_{i,n} - X_{i,n}(k).
\end{array}
\end{equation}
Both $X_{i,n}(k)$ and $\widetilde{X}_{i,n}(k)$ are continuous transformations of $X_{i,n}$ with Lipschitz constants 1. By Proposition \ref{proposition_transformation}, they inherit the dependence coefficient from $X_{i,n}$.

Since
\begin{align*}
    \expct\left|\sum_{i\in D_n} (X_{i,n} - \expct X_{i,n})\right| \leq &\expct\left|\sum_{i\in D_n} (X_{i,n}(k) - \expct X_{i,n}(k))\right| + \expct\left|\sum_{i\in D_n} (\widetilde{X}_{i,n}(k) - \expct \widetilde{X}_{i,n}(k))\right|\\
    \leq &\expct\left|\sum_{i\in D_n} (X_{i,n}(k) - \expct X_{i,n}(k))\right| + 2\sum_{i\in D_n}\expct|\widetilde{X}_{i,n}(k)|,
\end{align*} we have
\begin{align*}
     \norm{(|D_n|_c)^{-1}\sum_{i\in D_n} (X_{i,n} - \expct X_{i,n})}_1 \leq &\norm{(|D_n|_c)^{-1}\sum_{i\in D_n} (X_{i,n}(k) - \expct X_{i,n}(k))}_1\\
     &+ 2\sup_n\sup_{i \in D_n}\expct|\widetilde{X}_{i,n}(k)|
\end{align*}
Note that $\sup_n\sup_{i \in D_n}\expct|\widetilde{X}_{i,n}(k)| \leq \sup_n\sup_{i \in D_n}\expct\left[\left|X_{i,n}\right|\mathbf 1\{\left|X_{i,n}\right| \geq k\}\right]$ for any $k > 0$, then according to \eqref{unif_l1+_int}, it suffices to show that
\begin{equation}\label{proof_lln_eq1}
    \lim_{k\to\infty}\lim_{n\to\infty}\frac{1}{|D_n|_c}\norm{\sum_{i\in D_n} (X_{i,n}(k) - \expct X_{i,n}(k))}_1 = 0
\end{equation} in order to prove that \[\lim_{n\to\infty}\frac{1}{|D_n|_c}\norm{\sum_{i\in D_n} (X_{i,n} - \expct X_{i,n})}_1 = 0.\]

Let $\sigma^2_n(k) = \var\left[\sum_{i \in D_n} X_{i,n}(k)\right]$, then
\[\frac{1}{|D_n|_c}\norm{\sum_{i\in D_n} (X_{i,n}(k) - \expct X_{i,n}(k))}_1 \leq \frac{\sigma_n(k)}{|D_n|_c}\] by Lyapunov's inequality. Since $X_{i,n}(k)$ is a bounded function of $X_{i,n}$ with Lipschitz constant 1, Lemma A.1(iii) in \cite{jenish2009} and \eqref{cov_ineq_bounded} give \[\sigma^2_n(k) \leq C |D_n|_c\left[1+\sum_{r=1}^{\infty}r^{d-1}\bar\eta_{1,1}(r)\right].\] Recall from Assumption \ref{assumption_epsilon_lln} that $\bar\eta_{1,1}(r) = \bigO(r^{-\alpha})$ with $\alpha > d$. Therefore, for each $k > 0$,
\[
    \lim_{n\to\infty}\frac{\sigma_n(k)}{|D_n|_c}=0.
\]
This completes the proof.

\subsection{Proof of Theorem \ref{theorem_clt}}

Use the truncation decomposition in \eqref{truncation}. The variances of its two components are \[\sigma^2_n(k) = \var\left[\sum_{i \in D_n} X_{i,n}(k)\right], \qquad \widetilde{\sigma}^2_n(k) = \var\left[\sum_{i \in D_n} \widetilde{X}_{i,n}(k)\right].\]

\begin{claim}\label{clm_var_ineq}
    $|\sigma_n - \sigma_n(k)| \leq \widetilde{\sigma}_n(k).$
\end{claim}

\begin{proof}

Let \[S_n(k) = \sum_{i \in D_n}[X_{i,n}(k) - \expct X_{i,n}(k)], \qquad \widetilde{S}_n(k) = \sum_{i \in D_n}[\widetilde{X}_{i,n}(k) - \expct \widetilde{X}_{i,n}(k)].\] Note that $S_n = S_n(k) + \widetilde{S}_n(k)$, $\sigma_n = \|S_n\|_2$, $\sigma_n(k) = \|S_n(k)\|_2$ and $\widetilde{\sigma}_n(k) = \|\widetilde{S}_n(k)\|_2$, then the inequality could be derived according to Minkowski's inequality.

\end{proof}

Recalling from \eqref{unif_l2+_bound} that $\|X\|_m := \sup_n\sup_{i \in D_n}\|X_{i,n}\|_m < \infty$ for some $m > 2$, then for each $k > 0$, \[\|X(k)\|_m := \sup_n\sup_{i \in D_n}\|X_{i,n}(k)\|_m \leq \|X\|_m,\] and \[\|\widetilde{X}(k)\|_m := \sup_n\sup_{i \in D_n}\|\widetilde{X}_{i,n}(k)\|_m \leq \|X\|_m.\]

\begin{claim}\label{clm_var_bounds}
    There exist fixed constants $0 < C_* \leq C^* < \infty$ and $n_0\in\bN$ such that \[C_*|D_n|_c \leq \sigma^2_n \leq C^*|D_n|_c\] for all $n \geq n_0$.
\end{claim}

\begin{proof}

Assumption \ref{assumption_variance} implies that there exist $C_* > 0$ and $n_0\in\bN$ such that $C_*|D_n|_c \leq \sigma^2_n$ for all $n \geq n_0$, which proves the lower bound.

For the upper bound, according to the covariance inequalities \eqref{cov_ineq} derived in Lemma \ref{lma_cov_ineq} with $p = m$:
\begin{align*}
    \sigma_n^2 \leq &\sum_{i\in D_n}\expct X_{i,n}^2 + \sum_{\substack{i,j \in D_n\\ i \neq j}}|\cov(X_{i,n}, X_{j,n})|\\
    \leq &\|X\|_m^2|D_n|_c + C \|X\|_m^{\frac{m}{m-1}}\sum_{\substack{i,j \in D_n\\ i \neq j}}[\bar\eta_{1,1}(\rho(i,j))]^{\frac{m-2}{m-1}}\\
    \leq &\|X\|_m^2|D_n|_c + C \|X\|_m^{\frac{m}{m-1}}\sum_{i \in D_n}\sum_{h = 1}^{\infty}\sum_{\substack{j \in D_n\\ \rho(i,j)\in [h, h+1)}} [\bar\eta_{1,1}(\rho(i,j))]^{\frac{m-2}{m-1}}.
\end{align*} Lemma A.1 (iii) in \cite{jenish2009} gives \[\sup_{i \in D}\left|\{j \in D: \rho(i,j) \in [h, h+1)\}\right|_c \leq C h^{d-1}\] for $h \geq 1$. Therefore, there exists a constant $C^* > 0$ such that
\begin{align*}
    \sigma_n^2 \leq &\left\{\|X\|_m^2 + C \|X\|_m^{\frac{m}{m-1}}\sum_{h=1}^{\infty}h^{d-1}[\bar\eta_{1,1}(h)]^{\frac{m-2}{m-1}}\right\}|D_n|_c\\
    =: &C^*|D_n|_c,
\end{align*} where the last equality follows from Assumption \ref{assumption_eta_clt}(a).

\end{proof}

Observe that $\widetilde{X}_{i,n}(k)$ is a continuous function of $X_{i,n}$ with Lipschitz constant 1, therefore $\widetilde{X}_{i,n}(k)$ inherits the dependence coefficient from $X_{i,n}$ according to Proposition \ref{proposition_transformation}. For each $k > 0$,
\begin{align*}
    \widetilde{\sigma}^2_n(k) \leq &\sum_{i,j \in D_n}|\cov(X_{i,n} - X_{i,n}(k), X_{j,n} - X_{j,n}(k))|\\
    \leq &\sum_{\substack{i,j \in D_n\\ \rho(i,j) \leq r}}\left[|\cov(X_{i,n}, X_{j,n} - X_{j,n}(k))| + |\cov(X_{i,n}(k), X_{j,n} - X_{j,n}(k))|\right]\\
    &+\sum_{\substack{i,j \in D_n\\ \rho(i,j) > r}}\left[|\cov(X_{i,n}, X_{j,n} - X_{j,n}(k))| + |\cov(X_{i,n}(k), X_{j,n} - X_{j,n}(k))|\right]\\
    \leq &C r^d\|X\|_m^m k^{2-m}|D_n|_c + \left\{C\|X\|_m^{\frac{m}{m-1}}\sum_{h=r+1}^{\infty}h^{d-1}[\bar\eta_{1,1}(h)]^{\frac{m-2}{m-1}}\right\}|D_n|_c.
\end{align*} The last inequality follows from similar arguments in the proof of Lemma \ref{cov_ineq}, combined with Lemma A.1(ii), (iii) in \cite{jenish2009}. Let $r = k^{\delta}$, where $\delta \in (0, \frac{m-2}{d})$. Together with the lower bound for $\sigma^2_n$ in Claim \ref{clm_var_bounds}, this gives
\begin{equation}\label{proof_clt_eq4}
    \lim_{k \to \infty}\sup_{n \geq n_0}\frac{\widetilde{\sigma}^2_n(k)}{\sigma^2_n} = 0.
\end{equation} Combining Claim \ref{clm_var_ineq} with \eqref{proof_clt_eq4} we get
\begin{equation}\label{proof_clt_eq5}
    \lim_{k \to \infty}\sup_{n \geq n_0}\left|1 - \frac{\sigma_n(k)}{\sigma_n}\right| \leq \lim_{k \to \infty}\sup_{n \geq n_0}\frac{\widetilde{\sigma}_n(k)}{\sigma_n} = 0.
\end{equation}

On the other hand, note that $X_{i,n}(k)$ is a bounded function of $X_{i,n}$ with Lipschitz constant 1. By \eqref{cov_ineq_bounded} we have
\begin{equation}
    \sigma^2_n(k) \leq C |D_n|_c\left[1+\sum_{r=1}^{\infty}r^{d-1}\bar\eta_{1,1}(r)\right]
\end{equation} for each $k > 0$. Together with the lower bound for $\sigma^2_n$, this yields
\begin{equation}\label{proof_clt_eq6}
    \frac{\sigma_n(k)}{\sigma_n} \leq C < \infty
\end{equation} for every $k > 0$ and all $n \geq n_0$. This result and \eqref{proof_clt_eq5} play a key role in the following argument.

For the next step, we will adopt Lemma \ref{lma_brockwell} to reduce the problem of proving the CLT for the triangular array of random fields $X_{i,n}$ to that of proving the CLT for the bounded triangular array of random fields $X_{i,n}(k)$.
\begin{claim}
    We have
    \begin{equation}\label{result}
        \sigma^{-1}_n\sum_{i \in D_n}X_{i,n} \convd \mathcal N(0,1)
    \end{equation} if
    \begin{equation}\label{result_bounded}
        \sigma^{-1}_n(k)\sum_{i \in D_n}\left[X_{i,n}(k) - \expct X_{i,n}(k)\right] \convd \mathcal N(0,1)
    \end{equation} for each $k \in \bN$.
\end{claim}

\begin{proof}

Let $Z_n = \sigma^{-1}_n\sum_{i \in D_n}X_{i,n}$, and $V_{n,k} = \sigma^{-1}_n\sum_{i \in D_n}\left[X_{i,n}(k) - \expct X_{i,n}(k)\right]$. Let $\mu_n$ and $\nu$ be the probability measures of $Z_n$ and $V$ respectively. If $Z_n$ does not converge to $V$ in distribution, then the L\'{e}vy-Prokhorov metric $d(\mu_n, \nu)$ does not converge to 0 as $n \to \infty$, i.e. for any $\delta > 0$, there always exist sub-indices $(n_m)_{m \in \bN}$ such that $d(\mu_{n_m}, \nu) > \delta$ for all $n_m$. Next, we will find a sub-sequence of $(Z_{n_m})$ such that it converges to $V$, contradicting with that $d(\mu_{n_m}, \nu) > \delta$ for all $n_m$. By \eqref{proof_clt_eq6}, $\frac{\sigma_n(k)}{\sigma_n} \leq C$ for each $k \in \bN$ and all $n \geq n_0$. After discarding finitely many terms, assume that $n_m \geq n_0$. By the Bolzano--Weierstrass theorem we have:
\begin{itemize}
    \item For $k = 1$, there exists sub-sub-indices $(n_{m(l_1)})_{l_1 \in \bN}$ such that \[\lim_{l_1 \to \infty}\frac{\sigma_{n_{m(l_1)}}(1)}{\sigma_{n_{m(l_1)}}} = \alpha(1);\]
    \item For $k = 2$, there exists sub-sub-sub-indices $(n_{m(l_1(l_2))})_{l_2 \in \bN}$ such that \[\lim_{l_2 \to \infty}\frac{\sigma_{n_{m(l_1(l_2))}}(2)}{\sigma_{n_{m(l_1(l_2))}}} = \alpha(2);\] \[\cdots\]
\end{itemize} Now we could find a sub-sequence $(n^*_r)$ of $(n_m)$ by letting $n^*_1 = n_{m(1)}$, $n^*_2 = n_{m(l_1(2))}$, $\cdots$ such that \[\lim_{r\to\infty}\frac{\sigma_{n^*_r}(k)}{\sigma_{n^*_r}} = \alpha(k)\] for each $k \in \bN$.

Observe that \[V_{n^*_r,k} = \frac{\sigma_{n^*_r}(k)}{\sigma_{n^*_r}}\left\{\sigma^{-1}_{n^*_r}(k)\sum_{i \in D_{n^*_r}}\left[X_{i,n^*_r}(k) - \expct X_{i,n^*_r}(k)\right]\right\}.\] If \eqref{result_bounded} holds, then the first condition in Lemma \ref{lma_brockwell} is satisfied since \[V_{n^*_r,k} \convd V_k \sim \mathcal N(0, \alpha^2(k))\] as $r \to \infty.$ Recalling from \eqref{proof_clt_eq5}, \[\lim_{k\to\infty}|\alpha(k) - 1| \leq \lim_{k\to\infty}\lim_{r\to\infty}\left|\alpha(k) - \frac{\sigma_{n^*_r}(k)}{\sigma_{n^*_r}}\right| + \lim_{k\to\infty}\sup_{n \geq n_0}\left|\frac{\sigma_n(k)}{\sigma_n} - 1\right| = 0,\] hence the second condition in Lemma \ref{lma_brockwell} is also verified.

Using Markov's inequality, \[\prob(|Z_n - V_{n,k}| > \delta) = \prob\left(\left|\sigma^{-1}_n\sum_{i \in D_n}(\widetilde{X}_{i,n}(k) - \expct \widetilde{X}_{i,n}(k))\right| > \delta\right) \leq \frac{\widetilde{\sigma}^2_n(k)}{\delta^2\sigma^2_n},\] for any $\delta > 0$. Hence condition 3 in Lemma \ref{lma_brockwell} holds for $Z_n$ and $V_{n,k}$ because of \eqref{proof_clt_eq4}, obviously it also holds for the sub-sequences $Z_{n^*_r}$ and $V_{n^*_r,k}$.

Applying Lemma \ref{lma_brockwell} on the sub-sequences $Z_{n^*_r}$ and $V_{n^*_r,k}$, we have $Z_{n^*_r} \convd V$ as $r \to \infty$. Since $(n^*_r)$ is a sub-sequence of $(n_m)$, $Z_{n^*_r} \convd V$ contradicts with former assumption that $Z_n$ does not converge weakly to $V$.

\end{proof}

Now we consider the case when $(X_{i,n})$ are bounded as $\sup_n\sup_{i\in D_n}|X_{i,n}| \leq C_X$. Let $(d_n)_{n \geq 1}$ be a sequence such that $\lim_{n \to \infty}d_n = \infty$, $\lim_{n \to \infty}\frac{d_n^d}{(|D_n|_c)^{1/2}} = 0$, and \[\lim_{n\to\infty}\bar\eta_{\infty,1}(d_n)(|D_n|_c) = 0\] for $\eta$-coefficients. According to Assumption \ref{assumption_eta_clt}(b), we could set $d_n = (|D_n|_c)^q$ with $q \in (\frac{1}{\beta}, \frac{1}{2d})$ since $\beta > 2d$.

Define \[a_n = \sum_{\substack{i, j \in D_n\\ \rho(i,j) \leq d_n}}\cov(X_{i,n}, X_{j,n}).\] Recalling from the covariance inequality for bounded random variables \eqref{cov_ineq_bounded}, there exists constant $C > 0$ such that
\begin{equation}\label{proof_clt_eq7}
\begin{array}{llll}
    |\sigma_n^2 - a_n| &\leq \sum_{\substack{i, j \in D_n\\ \rho(i,j) > d_n}}|\cov(X_{i,n}, X_{j,n})|\\
    &\leq C\sum_{i \in D_n}\sum_{r = d_n}^{\infty}\sum_{\substack{j \in D_n\\ \rho(i,j)\in [r, r+1)}}\bar\eta_{1,1}(r)\\
    &\leq C |D_n|_c\sum_{r=d_n}^{\infty}r^{d-1}\bar\eta_{1,1}(r)\\
    &= o(|D_n|_c).
\end{array}
\end{equation} Then we have \[0 < \liminf_{n\to\infty}(|D_n|_c)^{-1}\sigma_n^2 \leq \liminf_{n\to\infty}(|D_n|_c)^{-1}a_n + \liminf_{n\to\infty}(|D_n|_c)^{-1}o(|D_n|_c).\] Through arguments similar to those in the proof of Claim \ref{clm_var_bounds}, we have $\sup_{n \geq n_0}a_n = \bigO(|D_n|_c)$. Consequently, $\sigma^2_n = a_n + o(|D_n|_c) = a_n[1 + o(1)]$ for sufficiently large $n$. Define \[\bar{S}_n = a_n^{-1/2}\sum_{i \in D_n}X_{i,n} = \frac{\sigma_n}{a_n^{1/2}}\sigma_n^{-1}\sum_{i \in D_n}X_{i,n},\] then it remains for us to show following convergence, which could be verified using Lemma \ref{lma_stein}.

\begin{claim}\label{clm_result_block}
    $\bar{S}_n \convd \mathcal N(0,1)$ as $n \to \infty$.
\end{claim}

\begin{proof}

The first condition in Lemma \ref{lma_stein} is satisfied since $a_n = \bigO(|D_n|_c)$ and $\sigma_n^2 = \bigO(|D_n|_c)$ for sufficiently large $n$. Then it suffices to verify the second condition, i.e.
\begin{equation}\label{proof_clt_eq8}
    \lim_{n\to\infty}\expct\left[(\img\lambda - \bar{S}_n)e^{\img\lambda\bar{S}_n}\right] = 0
\end{equation} for all $\lambda\in\bR$.

Let \[S_{i,n} = \sum_{\substack{j\in D_n\\\rho(i,j)\leq d_n}} X_{j,n}, \qquad \bar{S}_{i,n} = a_n^{-1/2}S_{i,n}.\] Then we can make decomposition as follows
\[(\img\lambda - \bar{S}_n)e^{\img\lambda\bar{S}_n} = T_{1,n} + T_{2,n} + T_{3,n},\] where
\begin{align*}
    &T_{1,n} = \img\lambda e^{\img\lambda\bar{S}_n}\left(1 - a_n^{-1}\sum_{i \in D_n}X_{i,n}S_{i,n}\right),\\
    &T_{2,n} = a_n^{-1/2}e^{\img\lambda\bar{S}_n}\sum_{i \in D_n}X_{i,n}\left(e^{-\img\lambda\bar{S}_{i,n}} + \img\lambda\bar{S}_{i,n} - 1\right),\\
    &T_{3,n} = -a_n^{-1/2}\sum_{i \in D_n}X_{i,n}e^{\img\lambda(\bar{S}_n - \bar{S}_{i,n})}.
\end{align*} For the next step, we will prove that $\lim_{n\to\infty}\expct|T_{k,n}| = 0$ for each $k = 1, 2, 3$.

We firstly consider the term $T_{1,n}$. Note that $\sum_{i \in D_n}\expct(X_{i,n}S_{i,n}) = a_n$, then for sufficiently large $n$ we have
\begin{align*}
    \expct|T_{1,n}|^2 = &\lambda^2\left[1 - 2a_n^{-1}\sum_{i \in D_n}\expct(X_{i,n}S_{i,n}) + a_n^{-2}\expct\left(\sum_{i \in D_n}X_{i,n}S_{i,n}\right)^2\right]\\
    = &\lambda^2\left[1 - 2a_n^{-1}a_n + a_n^{-2}\var\left(\sum_{i \in D_n}X_{i,n}S_{i,n}\right) + a_n^{-2}a_n^2\right]\\
    = &\lambda^2 a_n^{-2}\var\left(\sum_{\substack{i,j\in D_n \\ \rho(i,j) \leq d_n}}X_{i,n}X_{j,n}\right)\\
    = &\lambda^2 a_n^{-2}\sum_{\substack{i,j,k,l\in D_n \\ \rho(i,j) \leq d_n \\ \rho(k,l) \leq d_n}}\cov(X_{i,n}X_{j,n}, X_{k,n}X_{l,n})\\
    \leq &C_{\lambda}|D_n|_c^{-2}\sum_{\substack{i,j,k,l\in D_n \\ \rho(i,j) \leq d_n \\ \rho(k,l) \leq d_n \\ \rho(i,k) > 3d_n}}|\cov(X_{i,n}X_{j,n}, X_{k,n}X_{l,n})|\\
    &+ C_{\lambda}|D_n|_c^{-2}\sum_{\substack{i,j,k,l\in D_n \\ \rho(i,j) \leq d_n \\ \rho(k,l) \leq d_n \\ \rho(i,k) \leq 3d_n}}|\cov(X_{i,n}X_{j,n}, X_{k,n}X_{l,n})|,
\end{align*} for some $0 < C_{\lambda} < \infty$.

Define function $f_u: \bR^u \mapsto \bR$ as 
\begin{equation}\label{proof_clt_f_u}
    f(x_1, \cdots, x_u) = -C_X^u \vee x_1 \cdots x_u \wedge C_X^u,
\end{equation} then $f_u$ is a bounded Lipschitz function. Recalling from \eqref{cov_ineq_bounded} we have
\begin{align*}
    &|\cov(X_{i,n}, X_{j,n})| = |\cov(f_1(X_{i,n}), f_1(X_{j,n}))| \leq C \bar\eta_{1,1}(\rho(i,j)),\\
    &|\cov(X_{i,n}X_{j,n}, X_{k,n}X_{l,n})| = |\cov(f_2(X_{i,n},X_{j,n}), f_2(X_{k,n},X_{l,n}))| \leq C \bar\eta_{2,2}(\rho(\{i,j\},\{k,l\})),\\
    &|\cov(X_{i,n}, X_{j,n}X_{k,n}X_{l,n})| = |\cov(f_1(X_{i,n}), f_3(X_{j,n},X_{k,n},X_{l,n}))| \leq C \bar\eta_{1,3}(\rho(i,\{j,k,l\})).
\end{align*}

When $\rho(i,k) > 3d_n$, we have $\rho(\{i,j\}, \{k,l\}) > \rho(i,k) - 2d_n$. Let \[N_i(r) = \left|\left\{(j,k,l): \rho(i,j) \leq d_n, \rho(k,l) \leq d_n, 3d_n < r \leq \rho(i,k) < r + 1\right\}\right|_c.\] Then Lemma A.1(ii), (iv) in \cite{jenish2009} gives $\sup_{i\in\bR^d}N_i(r) \leq C d_n^{2d}r^{d-1}$. Hence, for each $i \in D_n$,
\begin{align*}
    &\sum_{\substack{j,k,l\in D_n \\ \rho(i,j) \leq d_n \\ \rho(k,l) \leq d_n \\ \rho(i,k) > 3d_n}}|\cov(X_{i,n}X_{j,n}, X_{k,n}X_{l,n})| \\
    \leq &C \sum_{r = 3d_n}^{\infty}\sup_{i\in\bR^d}N_i(r)\bar\eta_{2,2}(r - 2d_n)\\
    \leq &C d_n^{2d}\sum_{r = 3d_n}^{\infty}r^{d-1}\bar\eta_{2,2}(r - 2d_n)\\
    \leq &C d_n^{2d}\sum_{r = d_n}^{\infty}r^{d-1}\bar\eta_{2,2}(r).
\end{align*} Therefore,
\begin{equation}\label{proof_clt_eq9}
    \sum_{\substack{i,j,k,l\in D_n \\ \rho(i,j) \leq d_n \\ \rho(k,l) \leq d_n \\ \rho(i,k) > 3d_n}}|\cov(X_{i,n}X_{j,n}, X_{k,n}X_{l,n})| \leq C |D_n|_c d_n^{2d}.
\end{equation}

When $\rho(i,k) \leq 3d_n$, let $V_i(r)$ be the ball centered at $i$ with radius $r$. Then $V_i(4d_n)$ contains every $(j,k,l)$ satisfying $\rho(i,j) \leq d_n$, $\rho(k,l) \leq d_n$, and $\rho(i,k) \leq 3d_n$. Let \[M_i(r) = \left|\{(j,k,l): j,k,l \in V_i(4d_n),\ r \leq \rho(i, \{j,k,l\}) < r + 1\}\right|_c.\] Lemma A.1(ii), (v) in \cite{jenish2009} gives $\sup_{i\in\bR^d}M_i(r) \leq C d_n^{2d}r^{d-1}$ for $r\geq1$. Hence, for each $i \in D_n$,
\begin{align*}
    &\sum_{\substack{j,k,l\in D_n \\ \rho(i,j) \leq d_n \\ \rho(k,l) \leq d_n \\ \rho(i,k) \leq 3d_n}}|\cov(X_{i,n}X_{j,n}, X_{k,n}X_{l,n})| \\
    \leq &\sum_{j,k,l\in V_i(4d_n)}|\cov(X_{i,n}X_{j,n}, X_{k,n}X_{l,n})|\\
    \leq &\sum_{j,k,l\in V_i(4d_n)}\left[|\expct(X_{i,n}X_{j,n}X_{k,n}X_{l,n})| + |\expct(X_{i,n}X_{j,n})||\expct(X_{k,n}X_{l,n})|\right]\\
    \leq &C \sum_{j,k,l\in V_i(4d_n)}\left[\bar\eta_{1,3}(\rho(i,\{j,k,l\})) + \bar\eta_{1,1}(\rho(i,\{j,k,l\}))\right]\\
    \leq &C d_n^{2d} + C\sum_{r=1}^{4d_n}M_i(r)\left[\bar\eta_{1,3}(r)+\bar\eta_{1,1}(r)\right]\\
    \leq &C d_n^{2d}\left\{1+\sum_{r=1}^{4d_n}r^{d-1}\left[\bar\eta_{1,3}(r)+\bar\eta_{1,1}(r)\right]\right\}\\
    \leq &C d_n^{2d}.
\end{align*} By Assumption \ref{assumption_eta_clt}(a) we have
\begin{equation}\label{proof_clt_eq10}
    \sum_{\substack{i,j,k,l\in D_n \\ \rho(i,j) \leq d_n \\ \rho(k,l) \leq d_n \\ \rho(i,k) \leq 3d_n}}|\cov(X_{i,n}X_{j,n}, X_{k,n}X_{l,n})| \leq C |D_n|_c d_n^{2d}.
\end{equation}

Note that $\lim_{n\to\infty}\frac{d_n^{2d}}{|D_n|_c} = 0$. Thus, \eqref{proof_clt_eq9} and \eqref{proof_clt_eq10} imply that \[\expct|T_{1,n}|^2 \leq C C_{\lambda}\frac{d_n^{2d}}{|D_n|_c} \to 0\] as $n\to\infty$.

Now we consider the second term, for sufficiently large $n$ we have
\begin{align*}
    |T_{2,n}| = &|a_n^{-1/2}||\sum_{i \in D_n}X_{i,n}\left(e^{-\img\lambda\bar{S}_{i,n}} + \img\lambda\bar{S}_{i,n} - 1\right)|\\
    \leq &C C_X(|D_n|_c)^{-1/2}\sum_{i\in D_n}|e^{-\img\lambda\bar{S}_{i,n}} + \img\lambda\bar{S}_{i,n} - 1|.
\end{align*} Note that
\begin{align*}
    |\bar{S}_{i,n}| \leq &a_n^{-1/2}\sum_{\substack{j\in D_n \\ \rho(i,j) \leq d_n}}|X_{j,n}|\\
    \leq &C C_X a_n^{-1/2}d_n^d\\
    = &\bigO((|D_n|_c)^{-1/2})d_n^d.
\end{align*} The second inequality follows from Lemma A.1(ii) in \cite{jenish2009}. Then $\lim_{n\to\infty}|\bar{S}_{i,n}| = 0$, hence $|\img\lambda\bar{S}_{i,n}| < 1/2$ for sufficiently large $n$. Since $|e^{-z} + z - 1| \leq |z|^2$ for complex number $|z| < 1/2$, $|e^{-\img\lambda\bar{S}_{i,n}} + \img\lambda\bar{S}_{i,n} - 1| \leq \lambda^2|\bar{S}_{i,n}|^2$ a.s. for sufficiently large $n$.

Now we have
\begin{align*}
    \expct|T_{2,n}| \leq &C C_X(|D_n|_c)^{-1/2}\sum_{i\in D_n}\lambda^2\expct|\bar{S}_{i,n}|^2\\
    \leq &C C_X(|D_n|_c)^{1/2}\lambda^2\sup_{i\in D_n}\expct|\bar{S}_{i,n}|^2\\
    \leq &C C_X(|D_n|_c)^{1/2}\lambda^2 a_n^{-1} \sup_{i\in D_n}\sum_{\substack{j,k\in D_n \\ \rho(i,j)\leq d_n \\ \rho(i,k)\leq d_n}}|\expct(X_{j,n}X_{k,n})|\\
    \leq &C (|D_n|_c)^{-1/2}\sup_{i\in D_n}\sum_{\substack{j\in D_n \\ \rho(i,j)\leq d_n}}\left[C_X^2+\sum_{r=1}^{2d_n}N_j(r)\bar\eta_{1,1}(r)\right]\\
    \leq &C (|D_n|_c)^{-1/2}d_n^d\left[C_X^2+\sum_{r=1}^{2d_n}r^{d-1}\bar\eta_{1,1}(r)\right]\\
    \leq &C (|D_n|_c)^{-1/2}d_n^d,
\end{align*} where $N_j(r) = \left|\{k: r\leq \rho(j,k) < r+1\}\right|_c$. The last two inequalities follow from Lemma A.1(iii) in \cite{jenish2009} and Assumption \ref{assumption_eta_clt}(a). Then $\lim_{n\to\infty}\expct|T_{2,n}| = 0$.

As for the third term, we want to prove that $\lim_{n\to\infty}|\expct T_{3,n}| = 0$. Firstly note that 
\begin{align*}
    |\expct T_{3,n}| = &\left|\expct\left[a_n^{-1/2}\sum_{i \in D_n}X_{i,n}e^{\img\lambda(\bar{S}_n - \bar{S}_{i,n})}\right]\right|\\
    \leq &C(|D_n|_c)^{-1/2}\sum_{i\in D_n}\left|\expct X_{i,n}e^{\img\lambda(\bar{S}_n - \bar{S}_{i,n})}\right|\\
    \leq &C(|D_n|_c)^{-1/2}\sum_{i\in D_n}\left(\left|\expct X_{i,n}\cos{\lambda(\bar{S}_n - \bar{S}_{i,n})}\right| + \left|\expct X_{i,n}\sin{\lambda(\bar{S}_n - \bar{S}_{i,n})}\right|\right).
\end{align*} Let $f^*(\fX_i) = \cos{\lambda(\bar{S}_n - \bar{S}_{i,n})}$ where $\fX_i = (X_{j,n})_{j\in D_n, \rho(i,j) > d_n}$. $f^*$ is bounded with Lipschitz constant $\lip(f^*) = |\lambda|a_n^{-1/2}$, with domain $\bR^u$ for some $u \leq |D_n|_c$. Another bounded Lipschitz function $f_1$ is defined in $\eqref{proof_clt_f_u}$. By \eqref{dependence_coef} we have:
\begin{align*}
    |\expct [X_{i,n}\cos{\lambda(\bar{S}_n - \bar{S}_{i,n})}]| = &|\cov[f^*(\fX_i), f_1(X_{i,n})]|\\
    \leq &(C_X|D_n|_c|\lambda|a_n^{-1/2} + 1)\bar\eta_{\infty,1}(d_n).
\end{align*} Same holds if $f^*(\fX_i) = \sin{\lambda(\bar{S}_n - \bar{S}_{i,n})}$. Since $\lim_{n\to\infty}\bar\eta_{\infty,1}(d_n)|D_n|_c = 0$, we have
\[|\expct T_{3,n}| \leq C(|D_n|_c)^{-1/2}\sum_{i\in D_n}(|D_n|_c)^{1/2}\bar\eta_{\infty,1}(d_n) \leq C|D_n|_c\bar\eta_{\infty,1}(d_n) \to 0\] as $n\to\infty$.

\end{proof}

The proof of Theorem \ref{theorem_clt} is completed as Claim \ref{clm_result_block} is verified.

\subsection{Proof of results in Section \ref{section_dynamic_network_estimation}}\label{proof_dynamic_network_estimation}

\subsubsection{Limit-theory conditions for the dynamic-network triangular array of random fields}

For an arbitrary $\vartheta=(a,\beta,\gamma)\in\Theta$, let
$Z_{i,t,n}=X_{i,t,n}-\expct X_{i,t,n}$, which preserves the $\eta$-coefficients of $X_{i,t,n}$. \eqref{example_dynamic_network_causal_solution} gives
\[
    |Z_{i,t,n}|\leq\frac{1}{1-a}
    \leq\frac{1}{1-\bar a}.
\]
Together with \eqref{example_dynamic_network_eta_bound}, this verifies the
moment and dependence conditions of Theorems \ref{theorem_lln} and
\ref{theorem_clt}, uniformly in $\vartheta$. It remains to verify Assumption \ref{assumption_variance} for the CLT. Write
\begin{equation}\label{proof_dynamic_network_linear_sum}
\begin{aligned}
    \sum_{(i,t)\in D_n}Z_{i,t,n}
    &=\sum_{(i,t)\in D_n}\sum_{k=0}^{\infty}a^k
    \Bigl[\varepsilon_{A_k^\gamma(i,t),t-k,n}
    -p_{A_k^\gamma(i,t),t-k,n}(\beta)\Bigr]\\
    &=\sum_{(j,s)\in\bZ^2}c_{j,s,n}
    \bigl[\varepsilon_{j,s,n}-p_{j,s,n}(\beta)\bigr],
\end{aligned}
\end{equation}
where the second equality groups together all terms involving the same
innovation, with
\[
    c_{j,s,n}
    =\sum_{(i,t)\in D_n}\sum_{k=0}^{\infty}a^k
    \mathbf 1\{
    \bigl(A_k^\gamma(i,t),t-k\bigr)=(j,s)
    \}.
\]

Fix $(j,s)\in\bZ^2$ and $k\geq0$. Any observation to which
$\varepsilon_{j,s,n}$ contributes at lag $k$ must have index $(i,s+k)$
with $A_k^\gamma(i,s+k)=j$. Since every $\pi_t^\gamma$ is a bijection, the
map $i\mapsto A_k^\gamma(i,s+k)$ is also a bijection. Hence at most one such
observation belongs to $D_n$, and the contribution to $c_{j,s,n}$ at lag
$k$ is at most $a^k$. Therefore
\[
    0\leq c_{j,s,n}\leq\sum_{k=0}^{\infty}a^k
    =\frac{1}{1-a}.
\]
For each fixed $(i,t)\in D_n$ and $k\geq0$, exactly one index $(j,s)$
satisfies $\bigl(A_k^\gamma(i,t),t-k\bigr)=(j,s)$. Thus
\begin{equation}\label{proof_dynamic_network_coefficient_bounds}
\begin{aligned}
    \sum_{(j,s)\in\bZ^2}c_{j,s,n}
    &=\sum_{(i,t)\in D_n}\sum_{k=0}^{\infty}a^k\\
    &=\frac{|D_n|_c}{1-a}.
\end{aligned}
\end{equation}
Moreover, if $(j,s)\in D_n$, the term corresponding to $(i,t)=(j,s)$
and $k=0$ equals one because $A_0^\gamma(j,s)=j$. Since all other terms
are nonnegative, $c_{j,s,n}\geq1$.

Independence and \eqref{example_dynamic_network_prob} now imply
\begin{align*}
\var\left(\sum_{(i,t)\in D_n}Z_{i,t,n}\right)
&\geq\sum_{(j,s)\in D_n}
    c_{j,s,n}^2p_{j,s,n}(\beta)
    [1-p_{j,s,n}(\beta)]\\
&\geq p_0(1-p_0)|D_n|_c,
\end{align*}
while
\begin{align*}
\var\left(\sum_{(i,t)\in D_n}Z_{i,t,n}\right)
&\leq\frac14\sum_{(j,s)\in\bZ^2}c_{j,s,n}^2\\
&\leq\frac{1}{4(1-a)}
    \sum_{(j,s)\in\bZ^2}c_{j,s,n}\\
&=\frac{|D_n|_c}{4(1-a)^2}
\leq\frac{|D_n|_c}{4(1-\bar a)^2}.
\end{align*}
Assumption \ref{assumption_variance} is now verified for the centered triangular array of random fields, uniformly in $\vartheta$.

\subsubsection{Proof of Proposition \ref{proposition_dynamic_network_score_wd}}

For a finite block $U_n\subset D_n$, define
\[
    \cN_\gamma(U_n)
    =U_n\cup
    \left\{
    \bigl(A_1^\gamma(i,t),t-1\bigr):(i,t)\in U_n
    \right\}.
\]
If $U_n,W_n\subset D_n$ satisfy $\rho(U_n,W_n)\geq r$, every point in
$\cN_\gamma(U_n)$ or $\cN_\gamma(W_n)$ is within distance one of the
corresponding original block. The triangle inequality therefore gives
\begin{equation}\label{proof_estimation_expanded_distance}
    \rho\bigl(\cN_\gamma(U_n),\cN_\gamma(W_n)\bigr)\geq r-2.
\end{equation}
Moreover, $|\cN_\gamma(U_n)|_c\leq2|U_n|_c$.

Let $|U_n|_c=u$ and $|W_n|_c=v$, and take $f\in\cF_u$ and
$g\in\cG_v$, applied to the corresponding blocks of local vectors $(V_{i,t,n})_{(i,t)\in U_n}$ in \eqref{estimation_dynamic_network_local_vector}.
When regarded as functions of the distinct random coordinates in
$\cN_\gamma(U_n)$ and $\cN_\gamma(W_n)$, their sup-norms do not increase
and their Lipschitz constants increase by at most a factor of two. The
deterministic covariates $z_{i,t,n}$ do not alter these bounds. Hence,
by \eqref{dependence_coef}, \eqref{proof_estimation_expanded_distance},
and \eqref{example_dynamic_network_eta_bound}, for $r\geq3$,
\begin{align}
&\left|\cov\left[
f\bigl((V_{i,t,n})_{(i,t)\in U_n}\bigr),
g\bigl((V_{i,t,n})_{(i,t)\in W_n}\bigr)
\right]\right|
\notag\\
&\quad\leq4\Psi(f,g)
    \sup_{\vartheta\in\Theta}\bar\eta^\vartheta(r-2)
\notag\\
&\quad\leq\frac{8}{1-\bar a}
    \bar a^{\lceil(r-2)/2\rceil}\Psi(f,g).
\label{proof_estimation_local_eta_bound}
\end{align}
Dividing by $\Psi(f,g)$ and taking the relevant suprema gives the same
bound for the weak-dependence coefficient of the local-vector triangular array of random fields.

By \eqref{example_dynamic_network_causal_solution}, both $X_{i,t,n}$ and
$X_{A_1^\gamma(i,t),t-1,n}$ are bounded by $(1-\bar a)^{-1}$. The logistic function and
its first three derivatives are uniformly bounded when
$\sup_{i,t,n}\|z_{i,t,n}\|<\infty$ and $\xi\in\Xi$. Hence every scalar
statistic listed in Proposition \ref{proposition_dynamic_network_score_wd}
is uniformly bounded and Lipschitz as a function of $V_{i,t,n}$, with a
Lipschitz constant that is uniform in $\xi$, $c$, $k$, and $\ell$.
Proposition \ref{proposition_transformation} and \eqref{proof_estimation_local_eta_bound} now give \eqref{estimation_dynamic_network_statistic_eta}.

\subsubsection{Proof of Theorem \ref{theorem_dynamic_network_estimator}}

Since
\[
    \cF_{t-1,n}
    =\sigma\{\varepsilon_{j,s,n}:j\in\bZ,\ s\leq t-1\}.
\]
Then $X_{A_1^\gamma(i,t),t-1,n}$ is $\cF_{t-1,n}$-measurable and, at the true parameter $\xi_0$,
\begin{equation}\label{proof_estimation_true_residual}
\begin{aligned}
    e^0_{i,t,n}:=& e_{i,t,n}(\xi_0)
      =\varepsilon_{i,t,n}-p^0_{i,t,n},\\
    \expct(e^0_{i,t,n}\mid\cF_{t-1,n})=& 0,\\
    \expct[(e^0_{i,t,n})^2\mid\cF_{t-1,n}] =& p^0_{i,t,n}(1-p^0_{i,t,n}).
\end{aligned}
\end{equation}

We first prove consistency. Put
$\Delta_{i,t,n}(\xi)=\mu_{i,t,n}(\xi)-\mu_{i,t,n}(\xi_0)$.
Since $\Delta_{i,t,n}(\xi)$ is $\cF_{t-1,n}$-measurable,
\eqref{proof_estimation_true_residual} yields
\begin{equation}\label{proof_estimation_population_objective}
\begin{aligned}
    \expct Q_n(\xi)-\expct Q_n(\xi_0)
    &=\frac{1}{2|D_n|_c}\sum_{(i,t)\in D_n}
       \expct\Delta_{i,t,n}(\xi)^2.
\end{aligned}
\end{equation}
For every fixed $\xi\in\Xi$, Proposition
\ref{proposition_dynamic_network_score_wd} shows that the criterion
contributions are bounded and have exponentially decreasing
$\eta$-coefficients. Theorem \ref{theorem_lln} therefore gives
\begin{equation}\label{proof_estimation_pointwise_lln}
    Q_n(\xi)-\expct Q_n(\xi)\convLone0.
\end{equation}
The uniformly bounded first derivatives of the criterion contribution $\frac12e_{i,t,n}(\xi)^2$
give a constant $C$ such that, uniformly in $(i,t,n)$,
\[
    \left|\frac12e_{i,t,n}(\xi)^2
          -\frac12e_{i,t,n}(\widetilde\xi)^2\right|
    \leq C\|\xi-\widetilde\xi\|.
\]
Cover the compact set $\Xi$ by a finite $\delta$-net. Applying
\eqref{proof_estimation_pointwise_lln} at its finitely many points and
then letting $\delta\downarrow0$ proves
\begin{equation}\label{proof_estimation_uniform_lln}
    \sup_{\xi\in\Xi}|Q_n(\xi)-\expct Q_n(\xi)|\convp0.
\end{equation}
Combining \eqref{proof_estimation_population_objective},
\eqref{proof_estimation_uniform_lln}, and the separation condition
\eqref{estimation_dynamic_network_identification} gives
$\widehat\xi_n\convp\xi_0$.

We next establish CLT for the score $\nabla_\xi\left\{\frac12e_{i,t,n}(\xi)^2\right\}$. At $\xi_0$,
\begin{equation}\label{proof_estimation_true_score}
\begin{aligned}
    \left.\nabla_\xi\!\left\{\frac12e_{i,t,n}(\xi)^2\right\}
    \right|_{\xi=\xi_0}
    &=-e^0_{i,t,n}\nabla_\xi\mu_{i,t,n}(\xi_0),\\
    \expct\left[
    \left.\nabla_\xi\!\left\{
    \frac12e_{i,t,n}(\xi)^2\right\}
    \right|_{\xi=\xi_0}\right]&=0.
\end{aligned}
\end{equation}
For distinct indices, the score contributions in
\eqref{proof_estimation_true_score} are uncorrelated. Without loss of
generality, let $s\leq t$. If $s<t$, then
\begin{align*}
&\expct\left[
e^0_{i,t,n}e^0_{j,s,n}
\nabla_\xi\mu_{i,t,n}(\xi_0)
\nabla_\xi\mu_{j,s,n}(\xi_0)'
\right]\\
&\quad=
\expct\left[
e^0_{j,s,n}
\nabla_\xi\mu_{i,t,n}(\xi_0)
\nabla_\xi\mu_{j,s,n}(\xi_0)'
\expct(e^0_{i,t,n}\mid\cF_{t-1,n})
\right]
=0.
\end{align*}
If $s=t$ and $i\neq j$, conditional independence of the current
innovations gives
\begin{align*}
&\expct\left[
e^0_{i,t,n}e^0_{j,t,n}
\nabla_\xi\mu_{i,t,n}(\xi_0)
\nabla_\xi\mu_{j,t,n}(\xi_0)'
\right]\\
&\quad=
\expct\left[
\nabla_\xi\mu_{i,t,n}(\xi_0)
\nabla_\xi\mu_{j,t,n}(\xi_0)'
\expct(e^0_{i,t,n}e^0_{j,t,n}\mid\cF_{t-1,n})
\right]
=0.
\end{align*}
Moreover, \eqref{proof_estimation_true_residual} implies
\[
\begin{aligned}
&\expct\left[
(e^0_{i,t,n})^2
\nabla_\xi\mu_{i,t,n}(\xi_0)
\nabla_\xi\mu_{i,t,n}(\xi_0)'
\right]\\
&\quad=
\expct\left[
p^0_{i,t,n}(1-p^0_{i,t,n})
\nabla_\xi\mu_{i,t,n}(\xi_0)
\nabla_\xi\mu_{i,t,n}(\xi_0)'
\right].
\end{aligned}
\]
Summing these variance terms gives
\begin{equation}\label{proof_estimation_score_variance}
    \frac{1}{|D_n|_c}\var\left(
       \sum_{(i,t)\in D_n}
       \left.\nabla_\xi\!\left\{\frac12e_{i,t,n}(\xi)^2\right\}
       \right|_{\xi=\xi_0}\right)
    =\Omega_n.
\end{equation}
The score is bounded, and Proposition
\ref{proposition_dynamic_network_score_wd} provides exponential decay
uniformly over all block sizes. Thus the moment and dependence conditions
of Corollary \ref{corollary_clt} are satisfied. The positive definiteness
of the limit $\Omega$ verifies its variance condition through
\eqref{proof_estimation_score_variance}, since
\[
    \frac{1}{|D_n|_c}\lambda_{\min}\left[
      \var\left(\sum_{(i,t)\in D_n}
      \left.\nabla_\xi\!\left\{\frac12e_{i,t,n}(\xi)^2\right\}
      \right|_{\xi=\xi_0}\right)
    \right]
    =\lambda_{\min}(\Omega_n)\longrightarrow\lambda_{\min}(\Omega)>0.
\]
Corollary \ref{corollary_clt}
and Assumption \ref{assumption_dynamic_network_estimation}(c) imply
\begin{equation}\label{proof_estimation_score_clt}
    \frac{1}{\sqrt{|D_n|_c}}
       \sum_{(i,t)\in D_n}
       \left.\nabla_\xi\!\left\{\frac12e_{i,t,n}(\xi)^2\right\}
       \right|_{\xi=\xi_0}
    \convd \mathcal N(0,\Omega).
\end{equation}

For the Hessian $\nabla^2_{\xi\xi'}\!\left\{\frac12e_{i,t,n}(\xi)^2\right\}$, Proposition
\ref{proposition_dynamic_network_score_wd} and Theorem
\ref{theorem_lln}, applied entrywise, give a pointwise LLN. Uniformly
bounded third derivatives of the criterion make the Hessian uniformly
Lipschitz in $\xi$. The same finite-net argument used for
\eqref{proof_estimation_uniform_lln} therefore yields
\begin{equation}\label{proof_estimation_hessian_ulln}
    \sup_{\xi\in\Xi}
    \left\|
      \frac{1}{|D_n|_c}\sum_{(i,t)\in D_n}
      \left[
      \nabla^2_{\xi\xi'}\!\left\{\frac12e_{i,t,n}(\xi)^2\right\}
      -\expct\nabla^2_{\xi\xi'}\!\left\{\frac12e_{i,t,n}(\xi)^2\right\}
      \right]
    \right\|\convp0.
\end{equation}
Direct differentiation gives
\[
    \nabla^2_{\xi\xi'}\!\left\{\frac12e_{i,t,n}(\xi)^2\right\}
    =\nabla_\xi\mu_{i,t,n}(\xi)\nabla_\xi\mu_{i,t,n}(\xi)'
     -e_{i,t,n}(\xi)\nabla^2_{\xi\xi'}\mu_{i,t,n}(\xi).
\]
At the true value, \eqref{proof_estimation_true_residual} and the fact
that $\nabla^2_{\xi\xi'}\mu_{i,t,n}(\xi_0)$ is deterministic therefore give
\begin{equation}\label{proof_estimation_expected_hessian}
    \frac{1}{|D_n|_c}\sum_{(i,t)\in D_n}
      \expct\left[
      \left.\nabla^2_{\xi\xi'}\!\left\{
      \frac12e_{i,t,n}(\xi)^2\right\}
      \right|_{\xi=\xi_0}\right]
      =J_n\longrightarrow J.
\end{equation}
The uniform Lipschitz bound for the expected Hessian, consistency, and
\eqref{proof_estimation_hessian_ulln} imply
\begin{equation}\label{proof_estimation_integrated_hessian}
\begin{aligned}
    \cJ_n
    &:=\frac{1}{|D_n|_c}\sum_{(i,t)\in D_n}
       \int_0^1
       \left.\nabla^2_{\xi\xi'}\!\left\{
       \frac12e_{i,t,n}(\xi)^2\right\}
       \right|_{\xi=\xi_0+s(\widehat\xi_n-\xi_0)}\,ds
      \convp J.
\end{aligned}
\end{equation}
Because $\xi_0$ is an interior point, consistency places
$\widehat\xi_n$ in the interior of $\Xi$ with probability tending to
one. The first-order condition and the integral form of Taylor's theorem
then give
\[
    0=\frac{1}{|D_n|_c}\sum_{(i,t)\in D_n}
      \left.\nabla_\xi\!\left\{\frac12e_{i,t,n}(\xi)^2\right\}
      \right|_{\xi=\xi_0}
      +\cJ_n(\widehat\xi_n-\xi_0).
\]
Equations \eqref{proof_estimation_score_clt} and
\eqref{proof_estimation_integrated_hessian}, followed by Slutsky's
theorem, prove \eqref{estimation_dynamic_network_asymptotic_normality}.

Finally, \eqref{proof_estimation_hessian_ulln},
\eqref{proof_estimation_expected_hessian}, and consistency give
$\widehat J_n\convp J$. At the true value, Theorem \ref{theorem_lln}
and Proposition \ref{proposition_dynamic_network_score_wd}, applied entrywise
to
\[
    \left[
    \left.\nabla_\xi\!\left\{\frac12e_{i,t,n}(\xi)^2\right\}
    \right|_{\xi=\xi_0}\right]
    \left[
    \left.\nabla_\xi\!\left\{\frac12e_{i,t,n}(\xi)^2\right\}
    \right|_{\xi=\xi_0}\right]',
\]
show that its sample average minus its expectation converges to zero in
probability.
The average expectation equals $\Omega_n$ by
\eqref{proof_estimation_true_residual} and
\eqref{proof_estimation_true_score}. Uniform Lipschitz continuity in
$\xi$ and $\widehat\xi_n\convp\xi_0$ now yield
$\widehat\Omega_n\convp\Omega$. Since $J$ is positive definite,
$\widehat J_n$ is invertible with probability tending to one. Continuity
of matrix inversion completes
the proof of part (iii).

\newpage
\bibliographystyle{apalike}
\bibliography{references}

\end{document}